\documentclass[11pt,oneside]{article}

\usepackage{amscd}
\usepackage{amsfonts}
\usepackage[centertags]{amsmath}
\usepackage{amssymb}
\usepackage[italian,english]{babel}
\usepackage[all]{xy}
\usepackage{theorem}
\usepackage[dvips]{graphicx}

\usepackage[a4paper,hdivide={1.5cm,18.2cm, },dvips,bmargin=1.6cm]{geometry}

\newcommand{\Chi}{\ensuremath{\mathfrak X}}

\newcommand{\map}[3]{\mbox{${#1}\colon{#2}\to{#3}$}}
\newcommand{\dismap}[5]{
\[
\begin{array}{rcll}
#1: & #2 & \longrightarrow & #3\\
 & #4 & \longmapsto & #5
\end{array}
\]
}
\newcommand{\dismapnoname}[4]{
\[
\begin{array}{cll}
 #1 & \longrightarrow & #2\\
 #3 & \longmapsto & #4
\end{array}
\]
}

\newcommand{\R}{\ensuremath{\mathbb R}}
\newcommand{\C}{\ensuremath{\mathbb C}}
\newcommand{\Z}{\ensuremath{\mathbb Z}}

\newcommand{\Cp}[1]{\ensuremath{{\C\mathbb P}^{#1}}}

\newcommand{\ug}{\stackrel{\rm def}{=}}

\newcommand{\iso}{\ensuremath{\simeq}}
\newcommand{\tM}{\ensuremath{\tilde{M}}}
\newcommand{\tG}{\ensuremath{\tilde{G}}}
\newcommand{\tfund}{\ensuremath{\tilde{\Omega}}}
\newcommand{\tC}{\ensuremath{\tilde{\mathcal{C}}}}

\newcommand{\myref}[1]{(\ref{#1})}

\newcommand{\EndDim}{\ensuremath{\nopagebreak\hfill\blacksquare}}

\newcommand{\nul}{\ensuremath{\mathop{\text{Null\,}}}}

\newcommand{\gl}[1]{\ensuremath{\text{\upshape \rmfamily GL}(#1)}}
\newcommand{\pr}[1]{\ensuremath{\text{\upshape \rmfamily P}(#1)}}

\newcommand{\aut}[1]{\ensuremath{\text{\upshape \rmfamily Aut}(#1)}}

\newcommand{\im}[1]{\ensuremath{\text{\upshape \rmfamily Im}#1}}
\newcommand{\isom}[1]{\ensuremath{\text{\upshape \rmfamily Isom}(#1)}}
\newcommand{\cont}[1]{\ensuremath{\text{\upshape \rmfamily Cont}(#1)}}

\renewcommand{\setminus}{\smallsetminus}

\newtheorem{lem}{Lemma}[section]

\newtheorem{teo}[lem]{Theorem}
\newtheorem{cor}[lem]{Corollary}
\newtheorem{pro}[lem]{Proposition}
{\theorembodyfont{\rmfamily} \newtheorem{oss}[lem]{Remark}
                             \newtheorem{defi}[lem]{Definition}
                             \newtheorem{exa}[lem]{Example}
                            }

\newenvironment{D}[1][]{{\nopagebreak\noindent\em Proof#1: }}{\EndDim}

\newenvironment{acknowledgements}{{\bf Acknowledgements: }}{}
\newenvironment{address}{}{}

\newcommand{\lee}{\ensuremath{\omega}}
\newcommand{\fund}{\ensuremath{\Omega}}
\newcommand{\lie}{\ensuremath{\mathcal{L}}}
\newcommand{\liealg}[1]{\ensuremath{\mathfrak{#1}}}
\newcommand{\momenthat}{\ensuremath{\mu}}
\newcommand{\moment}{\ensuremath{\mu}}

\title{Locally conformal K\"ahler reduction}
\author{Rosa Gini\footnote{Partially supported by EURROMMAT.},
Liviu Ornea\footnote{Member of EDGE, Research
Training Network HRPN-CT-2000-00101, supported by the European Human
Potential Programme.}, Maurizio Parton$^{*,\dag}$
}

\begin{document}

\maketitle
\begin{abstract}
We define reduction of locally
conformal K\"ahler manifolds, considered as conformal Hermitian
manifolds, and we show its equivalence with an unpublished
construction given by Biquard and Gauduchon.
We show the compatibility between this reduction and
K\"ahler reduction of the universal cover. By a recent result of
Kamishima  and the second author, in the Vaisman case
(that is, when a metric in the conformal class has parallel Lee
form) if the manifold is compact
its universal cover comes equipped with the
structure of K\"ahler cone over a Sasaki compact manifold.
We show the compatibility between our reduction and Sasaki
reduction, hence describing
a subgroup of automorphisms whose action causes reduction to bear a Vaisman structure. Then we apply
this theory
to construct a wide class of Vaisman manifolds.
\end{abstract}



\section{Introduction}

Since $1974$ when the classical reduction procedure of S.~Lie was
formulated in modern terms by J.~Marsden and A.~Weinstein for symplectic
structures, this technique was extended to other various geometric
structures defined by a closed form. Extending the equivariant symplectic
reduction to K\"ahler manifolds was most natural: one only showed the
almost complex structure was also projectable. Generalizations to
hyperk\"ahler and quaternion K\"ahler geometry followed. The extension to
contact geometry is also natural and can be understood via the
symplectization of a contact manifold. In each case, the momentum map
is produced by a Lie group acting by specific automorphisms of the
structure.

A {\em locally conformal K\"ahler manifold} is a conformal Hermitian
manifold $(M,[g],J)$ such that for one (and hence for all) metric $g$
in the conformal class
the corresponding K\"ahler
form $\fund$ satisfies $d\fund=\lee\wedge\fund$, where $\lee$ is a closed
$1$-form.
This is equivalent to the existence of an atlas such that the
restriction of $g$ to any chart is conformal to a K\"ahler metric.

The $1$-form $\lee\in\Omega^1(M)$ was introduced by H.-C.~Lee in \cite{LeeKED},
and it is therefore called the {\em Lee form} of the
Hermitian structure $(g,J)$.

It was not obvious how to produce a quotient construction in conformal
geometry. The first published result we are aware of belongs to
S.~Haller and T.~Rybicki who proposed in \cite{HaRRLC} a reduction for
locally conformal symplectic structures. Their technique is
essentially local: they reduce the
local symplectic structures, then glue the local reduced structures.
But even earlier, since $1998$,
an unpublished paper by O.~Biquard and P.~Gauduchon proposed a
quotient construction for locally conformal K\"ahler manifolds
\cite{BiGCML}. Their construction relies heavily on the language and
techniques of conformal geometry as developed, for example, in
\cite{CaPEWG}. The key point is the fact that a locally conformal K\"ahler
structure can be seen as a closed $2$-form with values in a vector
bundle (of densities).

%The Haller-Rybicki construction, although
%completely independent and written in another language, is  equivalent
%with the locally conformal
%symplectic part of the Biquard-Gauduchon one.
%\dettagli{M: ho tolto il not surprising}

Our starting point was the paper \cite{HaRRLC}. Following the lines of
the K\"ahler reduction, we verified that the complex structure of a
locally conformal K\"ahler manifold can be projected to the quotient.
In section~\ref{momentmap} of this paper we construct the momentum map associated to
an action by locally conformal K\"ahler automorphisms, lying on the notion
of {\em twisted Hamiltonian action}\/ given
by I.~Vaisman in \cite{VaiLCS}. In section~\ref{theorem}
we extend Haller-Rybicki construction to the complex setting. Then,
in section~\ref{biqgau} we present,
rather in detail, due to its very restricted previous circulation, the
Biquard-Gauduchon construction. The main result of this section proves the
equivalence between the Biquard-Gauduchon reduction and ours.

The universal cover of a locally conformal K\"ahler manifold has a
natural  (global) homothetic K\"ahler structure. We exploit this
fact in section~\ref{trivialclass} in order to relate locally conformal K\"ahler reduction
to the K\"ahler reduction of its universal cover.

The study of locally conformal K\"ahler manifolds started in the
field of Hermitian manifolds. Most of the known examples of locally conformal K\"ahler metrics are
on compact manifolds and enjoy the additional property of having
parallel Lee form with respect to the Levi-Civita connection.
Locally conformal K\"ahler metrics with parallel Lee form were first
introduced and studied by I.~Vaisman in~\cite{VaiLCK,VaiGHM}, so we call {\em Vaisman metric}\/ a
locally conformal K\"ahler metric with this property.
 %under the name of {\em
%generalized Hopf manifolds}. We prefer to call them {\em Vaisman
%manifolds}.
Manifolds bearing a Vaisman metric show a rich geometry.
Such are
the Hopf surfaces $H_{\alpha, \beta}$ described in \cite{GaOLCK},
all
diffeomorphic with $S^1\times S^3$ (see also \cite{ParHSL}). More
generally, I.~Vaisman firstly showed that on the product $S^1\times
S^{2n+1}$ given as a quotient of $\C^n\setminus 0$ by the cyclic
infinite group spanned by $z\mapsto\alpha z$, where $z\in\C^n\setminus
0$ and $|\alpha|\neq 1$,
the projection of the metric $|z|^{-2}\sum dz^i\otimes d\bar z^i$ is locally conformal K\"ahler
with parallel Lee form $-|z|^{-2}\sum (z_id\bar z_i
+\bar z_idz_i)$.
The complete list
of compact complex locally conformal K\"ahler surfaces admitting parallel Lee form was given by
F.~Belgun in \cite{BelMSN} where it is also proved the existence of some
compact complex surfaces  which do not admit any locally conformal
K\"ahler metric.

The definition of Vaisman metric is {\em not}\/ invariant up to
conformal changes. A conformally equivalent notion of Vaisman manifold is still missing, but
a recent result by Kamishima and the second author in \cite{KaOGFC} provides one in the
compact case, generalizing
the one first proposed by Belgun in \cite{BelMSN} in the case of surfaces.
 We
develop this notion in section~\ref{vaisman} where we analyze
reduction in this case.

Vaisman geometry is  closely related with Sasaki
geometry. In this case the picture turns out to be the following. The category of ordinary locally
conformal K\"ahler manifolds can be seen as the image of the
category of pairs $(K,\Gamma)$ of homothetic
K\"ahler manifolds with a subgroup $\Gamma$ of homotheties acting freely
and properly discontinuously, with
morphisms given by homothetic K\"ahler morphisms equivariant by
the actions. What we prove in section~\ref{trivialclass} is that
under the functor associating to $(K,\Gamma)$ the locally
conformal K\"ahler manifold $K/\Gamma$ Hamiltonian actions go to twisted Hamiltonian actions, and vice versa,
see Theorem~\ref{compatibilitakaehler}. So
the images of subgroups
producing K\"ahler reduction actually are subgroups producing
locally conformal K\"ahler reduction (up to topological
conditions), and vice versa.
The same way the category of (compact) Vaisman
manifolds can be seen as the image of the category of pairs
$(W,\Gamma)$, with $W$ a Sasaki manifold and $\Gamma$ a subgroup
of (proper) homotheties of the K\"ahler cone $W\times\R$ acting
freely and properly discontinuously, with morphisms given by
Sasaki morphisms equivariant by the actions. The functor
associating to $(W,\Gamma)$ the Vaisman manifold
$(W\times\R)/\Gamma$
is
surjective on objects but not on morphisms: we call {\em Vaisman
morphisms}\/ the ones in the image. What we prove in
section~\ref{vaisman} is that, up to topological conditions,
subgroups of Sasaki automorphisms producing Sasaki reduction go to
subgroups producing Vaisman reduction, and vice versa. This is
particularly remarkable since, up to topological conditions, Sasaki reduction applies to {\em any}\/
subgroup of automorphisms, that is, the momentum map is {\em
always} defined. So we obtain that reduction by Vaisman
automorphisms is always defined (up to topological conditions) and
produces Vaisman manifolds.

This allows building a wide set of Vaisman manifolds, reduced by circle
actions on Hopf manifolds, in section~\ref{vai}.

\begin{acknowledgements}
We thank F.~Belgun, C.~Boyer, E.~Ferrand, K.~Galicki, P.~Gauduchon, H.~Pe\-der\-sen, T.~S.~Ra\-tiu,
A.~Swann
for many enlightening discussions during the elaboration of this paper.

The very beginning of this paper arises to the visit of the third
named author to the \'Ecole Polytechnique in Paris, and he wishes in particular to thank Paul
Gauduchon for giving him the opportunity to read his unpublished
paper.

This work was then initiated during the visit of the first and
third named
author to the Institute of Mathematics  {\em ``Simion Stoilow''} of Bucharest.
They both would like to thank
 Vasile Br\^{\i}nzanescu  and Paltin %{\em doctor humoris causa}\/
 Ionescu for warm hospitality.

Moreover, this work was ended during the visit of the third named
author in Odense, Denmark. He wishes to thank in particular Henrik Pedersen
and Andrew Swann, together with their families, for making his permanence so beautiful.
\end{acknowledgements}

\section{Locally conformal K\"ahler manifolds}

Let $(M,J)$ be any almost-complex $n$-manifold, $n\geq4$, let $g$
be a Hermitian metric on $(M,J)$.
Let $\fund$ be the
K\"ahler form  defined by $\fund(X,Y)\ug g(JX,Y)$.
The map \map{L}{\Omega^1(M)}{\Omega^3(M)}
given by the wedging with $\fund$ is injective, so that the $g$-orthogonal
splitting $\Omega^3(M)=\im{L}\oplus(\im{L})_0$ induces a well-defined
$\lee\in\Omega^1(M)$ given by the relation
$d\fund=\lee\wedge\fund+(d\fund)_0$.
The $1$-form $\lee\in\Omega^1(M)$ is called the {\em Lee form}\/ of the
almost-Hermitian structure $(g,J)$.

A relevant notion in this setting is that of {\em twisted
differential}. Given a $p$-form $\psi$ its twisted differential is
the $(p+1)$-form
\[
d^\lee\psi\ug d\psi
-\lee\wedge\psi.
\]
Remark that $d^\lee\circ d^\lee=0$ if and only if $d\lee=0$.

A Hermitian metric $g$ on a complex manifold $(M,J)$ is said to be {\em locally conformal
K\"ahler}\/ if $g$ is (locally) conformal to local K\"ahler metrics. In this case
the local forms $d \alpha_U$ coming from the local
conformal factors $e^{\alpha_{U}}$ paste to a
global form $\lee$ satisfying $d\fund=\lee\wedge\fund$. Vice versa
this
last equation together with $d\lee=0$ characterizes the locally
conformal K\"ahler metrics. In other words a Hermitian metric
is locally conformal K\"ahler if and only if
\begin{equation}\label{itb}
d^\lee\circ d^\lee=0 \qquad \text{and} \qquad d^\lee\fund=0.
\end{equation}

\begin{defi}
A conformal Hermitian manifold $(M,[g],J)$ of complex dimension bigger than
$1$
is said to be a {\em locally  conformal K\"ahler
manifold}\/ if one (and hence all of) the metrics in $[g]$ is
locally conformal K\"ahler.
\end{defi}

\begin{oss}
If, in particular, the Lee form of one
(and hence all) of the metrics in $[g]$ is exact, then the manifold
is said to be {\em globally conformal K\"ahler}. This is in fact
equivalent to requiring that in the conformal class there exists a
K\"ahler metric, that is, any metric in $[g]$ is globally conformal to a
K\"ahler metric. From~\cite{VaiLGC} it is known that for
compact manifolds possessing a K\"ahler structure forbids
existence of locally non-globally conformal K\"ahler structures,
so the two worlds are generally considered as disjoint.
 In
this paper, however, the two notions behave the same way, so we consider the
global case as a subclass of the local case.
\end{oss}

>From now on, let $(M,[g],J)$ be a locally conformal K\"ahler manifold.

Not unlike  the K\"ahler case,
locally conformal K\"ahler manifolds come equipped with a notable subset
of $\Chi(M)$: given a smooth function $f$
the {\em associated  Hamiltonian vector field} is the $\fund$-dual of
$df$, and
{\em Hamiltonian vector fields} are vector fields that admit such
a presentation.
But the notion that works for reduction, as shown in
\cite{HaRRLC}, is the one given in \cite{VaiLCS} obtained by twisting the classical.
Given $f$\/ its {\em associated  twisted Hamiltonian vector field\/} is
 the \fund-dual of $d^\omega f$. The subset of $\Chi(M)$ of {\em twisted Hamiltonian vector
fields}\/ is that of vector fields admitting such a presentation.

\begin{oss}\label{unica}
If $M$ is not globally conformal K\"ahler the function associating
to $f$ its twisted Hamiltonian vector field is injective. Indeed
$d^\lee f=0$ implies $\lee=d\log|f|$ on $f\neq0$, so either $f\equiv0$
or $\lee$ is exact.
\end{oss}

Define a {\em twisted Poisson bracket}\/ on
$C^\infty(M)$ by
\begin{equation}\label{poisson}
\{f_1,f_2\}\ug\fund(\sharp d^\lee f_1,\sharp d^\lee f_2)
\end{equation}
The
relation
\begin{equation}\label{lieid}
\{\{f_1,f_2\},f_3\}+\{\{f_2,f_3\},f_1\}+\{\{f_3,f_1\},f_2\}=d^\lee\fund(\sharp d^\lee f_1,\sharp d^\lee f_2,\sharp d^\lee
f_3)=0
\end{equation}
proves that this bracket turns $C^\infty(M)$ into a Lie
algebra. Remark that the first equality in \eqref{lieid} holds
generally on any almost-Hermitian manifold $(M,g,J)$
under the only assumption $d\lee=0$.

\begin{oss}\label{conformalhamiltonian} Remark that the notion
of Hamiltonian vector field is invariant up to conformal
change of the metric, even though the function (possibly, the functions)
associated to a Hamiltonian vector field
 changes by the conformal factor.  A
straightforward computation shows in fact that, if
$\fund'=e^\alpha\fund$ and $\lee'=\lee+d\alpha$ is the
corresponding Lee form, the following relations hold
\[
\begin{split}
d^\lee f&=e^{-\alpha}d^{\lee'}(e^\alpha f)\\
\sharp_\fund d^\lee f&=\sharp_{\fund'}d^{\lee'}(e^\alpha f)\\
\{e^\alpha f_1,e^\alpha f_2\}^{\fund'}&=e^\alpha \{f_1,f_2\}^{\fund}
\end{split}
\]
so that multiplication by $e^\alpha$ yields an isomorphism between
$(C^\infty(M),\{\ ,\ \}^\fund)$ and $(C^\infty(M),\{\ ,\
\}^{\fund'})$ commuting with the corresponding maps $\sharp_\fund
d^\lee$ and $\sharp_{\fund'}d^{\lee'}$ in the space of twisted Hamiltonian
vector fields.
 In particular if $M$ is { globally conformal K\"ahler},
 then the twisted Hamiltonian vector fields of $M$ coincide with the
ordinary Hamiltonian vector fields, since the Lee form of a
K\"ahler metric is $0$.
\end{oss}

\begin{defi}
Given two locally conformal K\"ahler manifolds $(M,[g],J)$ and
$(M',[g'],J)$ a smooth map $h$ from $M$ to $M'$ is a {\em locally conformal K\"ahler
morphism}\/ if $h^*J'=J$ and $[h^*g']=[g]$. We denote by $\aut{M,[g],J}$, or briefly by $\aut{M}$, the
group of locally conformal K\"ahler automorphisms of $(M,[g],J)$.
\end{defi}

The group \aut{M}\ is a Lie group, contained as a subgroup in the complex Lie
group of biholomorphisms of $(M,J)$. However, unlike the Riemannian
case,
 the Lie algebra of
\aut{M} is {\em not}\/ closed for the complex structure. This will be used
in the sequel.

\section{The locally conformal K\"ahler momentum map}\label{momentmap}

In this paper we consider (connected) Lie
subgroups $G$
of $\aut{M}$.

\begin{oss}
It follows from \cite{MPPCEW} that whenever a locally conformal
K\"ahler manifold
$M$ is compact, the
group \aut{M}\ coincides with the isometries
of the {\em Gauduchon metric}\/ in  the conformal class, that is, the one
whose Lee form is coclosed. Hence, in particular, \aut{M} is compact.
More generally if a
subgroup $G$ of \aut{M}\ is compact then by using the Haar integral one obtains a metric in
the conformal class such that $G$ is contained in the group of its
isometries.
So the case when $G$ is {\em not}\/ constituted by isometries
of a specific metric can only happen if both $M$ and $G$ are
non-compact.
\end{oss}

Throughout the paper we identify
fundamental vector fields with elements $X$ of the
Lie algebra $\liealg{g}$ of $G$, so that if $x\in M$ then $\liealg{g}(x)$
means $T_x(Gx)$.

%Recall that in the symplectic
%Marsden-Weinstein reduction, where symplectomorphisms are
%used, the fundamental vector fields $X\in\liealg{g}$  are infinitesimal
%transformations
%of $\fund$, that is, $\lie_X\fund=0$. In our case, to require that the group acts by
%conformal and complex maps implies that $\lie_X^\lee\fund=c_X\fund$, for
%locally constant functions $c_X$.

Imitating the
terminology established in \cite{McSIST}, we call the
action of $G$ {\em weakly twisted  Hamiltonian}\/ if the associated infinitesimal
action is of twisted Hamiltonian vector fields, that is, if
there exists a (linear) map
\map{\momenthat^\cdot}{\liealg{g}}{C^\infty(M)} such that
$\iota_X\fund=d^\lee\momenthat^X$ for fundamental vector fields
$X\in\liealg{g}$, and {\em twisted Hamiltonian}\, if $\momenthat$ can be chosen to be a Lie
algebra homomorphism with respect to the Poisson bracket
\eqref{poisson}. %Remark that, since $d^\lee\circ d^\lee=0$, the
%action can be weakly Hamiltonian only if $c_X=0$, for
%$X\in\liealg{g}$.
In this case we say that the Lie
algebra homomorphism $\momenthat$ is a {\em momentum map} for the action
of $G$, or, equivalently, with the same name and symbol we refer to the induced map \map{\moment}{M}{\liealg{g}^*}
given by $\langle\moment(x),X\rangle \ug\momenthat^X(x)$, for $X\in \liealg{g}$ and
carets denoting the evaluation.

\begin{oss}\label{conformalpoisson}
Note that the property of an action
of being twisted Hamiltonian is
a property of the conformal structure, even though the
Poisson structure on $C^\infty(M)$ is not conformally invariant, see
Remark~\ref{conformalhamiltonian}. If $g'=e^\alpha g$ then
$\moment^{\lee'}=e^\alpha\moment^{\lee}$. In particular the preimage
of $0$ is well-defined.
\end{oss}

\begin{oss}
Remark that $\moment$ is not equivariant for the standard coadjoint action on $\liealg{g}^*$.
It is known from \cite{HaRRLC} that by
 modifying the coadjoint action by means of  the conformal factors arising from $h^*g\sim g$
one can force $\moment$ to be equivariant.
\end{oss}

On $\moment^{-1}(0)$ the twisted differential of
the associated twisted Hamiltonian functions $\moment(\mathfrak{g})$ coincides with the ordinary
differential, since
$d_x^\lee\momenthat^X=d_x\momenthat^X-%\stackrel{=0}
{\momenthat^X(x)}\lee_x$
for $X\in\liealg{g}$, $x\in\moment^{-1}(0).$
Thus, if the action is
twisted Hamiltonian, then
the
functions in $\moment(\mathfrak{g})$ vanish on the whole orbit of
$x\in\moment^{-1}(0)$, since
for $x\in\moment^{-1}(0)$ and  $Y\in\liealg{g}(x)$
\[
d_x\momenthat^X(Y)
=d_x^\lee\momenthat^X(Y)
=\fund(\sharp d^\lee\momenthat^X,\sharp d^\lee\momenthat^Y)(x)=\{\momenthat^X,\momenthat^Y\}(x)
=\momenthat^{[X,Y]}(x)=
0,
\]
that is to say, $\moment^{-1}(0)$
is closed for the action of $G$.

Moreover, if $0$ is a regular value for $\moment$, then
$T(\moment^{-1}(0))^{\perp_\fund}=\liealg{g}$, since
for any $x\in\moment^{-1}(0)$, $X\in\liealg{g}, V\in
\Chi(\moment^{-1}(0))$ we have
\[
\fund(X,V)(x)=d_x^\lee\momenthat^X(V)=d_x\momenthat^X(V)=0.
\]
Thus we say
that $\moment^{-1}(0)$ is a {\em coisotropic} submanifold of $M$.

In the next section we show how to obtain a locally conformal
K\"ahler structure on $\moment^{-1}(0)/G$ under the additional
hypothesis of it being a manifold. But we remark here that,
due to the missing
equivariance of $\moment$,
a non-zero reduction  is not
straightforward.

\begin{oss}
We give a brief description of the existence and unicity
problem for momentum maps.
Suppose the action is weakly twisted Hamiltonian, and choose a linear map
\map{\momenthat^\cdot}{\liealg{g}}{C^\infty(M)}.  Denote by $N$ the kernel
of \map{d^\lee}{C^\infty(M)}{\Omega^1(M)}.
The
obstruction for $\momenthat$ to be a Lie algebra homomorphism is
given by the map
\map{\tau}{\liealg{g}\times\liealg{g}}{N}
sending $(X,Y)$ into $\{\momenthat^X,\momenthat^Y\}-\momenthat^{[X,Y]}$,
which can be shown to live in $H^2(\liealg{g},N)$, and this
cohomology class vanishes whenever the action is twisted Hamiltonian. If
this is the case, then momentum maps are
parameterized by $H^1(\liealg{g},N)$.
If $(M,[g],J)$ is non-globally conformal K\"ahler,
then $N=0$, see Remark~\ref{unica}.
%
%Together with the classical result of
%Vaisman (see \cite{VaiLGC}) that a locally conformal K\"ahler
%metric on a compact manifold is conformal K\"ahler if and only if the manifold
%admits a K\"ahler metric, the above discussion
%implies that
Then, in particular, a weakly Hamiltonian action on a compact non-K\"ahler locally
conformal K\"ahler manifold always admits a
unique momentum map.
\end{oss}

In the following we will often need a technical lemma we prove here
once and for all. If $g$ and $g'$ are tensors on the same manifold we write $g\sim
g'$ if they are conformal to each other.

\begin{lem}\label{incollare} Let $M$ be a manifold, let $\{U_i\}_{i\in\mathcal{I}}$
be a locally finite open covering. Let  $\{\rho_i\}$ be a partition of unity relative to
$\{U_i\}$. The following three facts hold.

\begin{enumerate}

\item Let $g$ and $g'$ be two tensors globally defined on $M$ and such that
for any $i$
\[
g|_{U_i}\sim g'|_{U_i};
\]
then $g$ and $g'$ are globally conformal.
\item Let $\{g_i\}$ be a collection of local tensors, where
$g_i$ is defined on $U_i$, such that whenever $U_i\cap U_{j}\neq\emptyset$
\[
g_i|_{U_i\cap U_{j}}\sim g_{j}|_{U_i\cap U_{j}};
\]
then the tensor $g\ug\sum_{i}\rho_i g_{i}$ is globally
defined on $M$ and $g|_{U_i}$ is
locally conformal to $g_i$.
%
%\item Let $\{g_i\}$ is a collecas in ii) and $g$ is
%a global tensor such that $g$ is locally conformal to one of the
%$g_i$'s, then the tensor global $g'$ obtained as in i) is
%globally conformal to $g$.

\item Let $\{g_i\}$ and $g$ be as  in ii). If  $g'$ is
a global tensor such that $g'|_{U_i}$ is locally conformal to
$g_i$, then $g$ and $g'$ are globally conformal.

\end{enumerate}
\end{lem}

\begin{D}
First prove {\em i)}. Let $e^{\alpha_i}$ be the conformal factor such that
\[
g|_{U_i}=e^{\alpha_i} g'|_{U_i};
\]
then recalling that $\sum_i\rho_i=1$ one obtains
\[
g=(\sum_i\rho_i e^{\alpha_i}) g'.
\]

Now turn to {\em ii)}.  For any $x \in M$ let $U_x$ be
a neighborhood of $x$ which is completely contained in any
$U_i$ that contains $x$, let $U_{i_x}$ be one of them and
$e^{\alpha_{x,i}}$ be the conformal factor between $g_{i_x}$ and
$g_{i}$, defined  on $U_{i_x}\cap U_{i}$ which contains $U_x$:
then the following holds
\[
g|_{U_x}=(\sum_{i}\rho_{i}e^{\alpha_{x,i}})g_{i_x}.
\]

Finally {\em i)} and {\em ii)} imply {\em iii)}.
\end{D}

\begin{oss}
Using a more sophisticated argument it is proved in \cite{HaRRLC}
that in case {\em ii)} one obtains $g|_{U_i}\sim
g_i$.
\end{oss}

\section{The reduction theorem}\label{theorem}

\begin{teo}\label{teonostro}
Let $(M,[g],J)$ be a locally conformal K\"ahler manifold.
Let $G$ be a Lie subgroup of $\aut{M}$ whose action is twisted Hamiltonian and is free and
proper on $\moment^{-1}(0)$, $0$ being a regular value for  the
momentum map
$\moment$.
Then
there exists a locally conformal K\"ahler structure $([\bar{g}],\bar{J})$ on
$\moment^{-1}(0)/G$, uniquely determined by the
condition $\pi^* \bar{g}\sim i^*g$,
where $i$ denotes the inclusion of $\moment^{-1}(0)$ into $M$ and
$\pi$ denotes the projection of $\moment^{-1}(0)$ onto its
quotient.
\end{teo}

\begin{D}
Since $\moment^{-1}(0)$ is coisotropic, and its isotropic leaves are the
orbits of $G$, the $[g]$-orthogonal
splitting $T_xM=E_x\oplus\liealg{g}(x)\oplus J\liealg{g}(x)$
holds, where $E_x$ is the $[g]$-orthogonal complement of $\liealg{g}(x)$ in
$T_x(\moment^{-1}(0))$. This shows that $E$ is a complex subbundle of
$TM$ and, since $J$ is constant along $\liealg{g}$, it induces an
almost complex structure $\bar{J}$ on $\moment^{-1}(0)/G$. This is proven to be
integrable the same way as in the K\"ahler case,
by computing the Nijenhuis tensor of $\bar J$ and recalling that
$\pi_*[V,W]=[\pi_*V,\pi_*W]$ for projectable
vector fields $V,W$.

Take an open cover $\mathcal{U}$ of $\moment^{-1}(0)/G$ that trivializes
the $G$-principal bundle $\map{\pi}{\moment^{-1}(0)}{\moment^{-1}(0)/G}$
and for each $U\in \mathcal{U}$ choose a local section
$s_U$ of $\pi$.

Fix an open set $U$. On its preimage we have two
horizontal distributions: the (global) already defined distribution $E$,
$[g]$-horizontal, and the tangent distribution $S_U$ to $s_U(U)$, translated
along the fibres by means of the action of $G$ to give a
distribution on the whole preimage of $U$. Remark that $S_U$
cannot be chosen to coincide with $E$ in general, since $S_U$ is obviously a
(local) foliation, whereas
$E$ is not
integrable in general.

Given a vector field $\bar{V}$ on $U$ denote by $V$ its
$[g]$-horizontal lifting. Then for any $\bar{V}$ the  vector fields $V$
and $J(V)$
are projectable and
 $\bar{J}(\bar{V})=\pi_*J(V)$. Moreover denote by $V+\nu_V$ the lifting of $\bar{V}$ tangent to
$S_U$, so that   $ds_U(\bar{V})=V+\nu_V$.\footnote{ To be  precise we
should write this expression in the form
$ds_U(\bar{V})=V\circ s_U+\nu_V\circ s_U$.  }
Remark that $\nu_V$ is a vertical vector field on
$\pi^{-1}(U)$, and that clearly $V+\nu_V$ is projectable itself: more
explicitly, for a generic $x\in\pi^{-1}(U)$,
$$
(V+\nu_V)_x=(h_x^{-1})_*d_{\pi(x)}s_U(\bar{V}_{\pi(x)})
$$
where by $h_x$ we denote the element of $G$ that takes
$x$ in $s_U(\pi(x))$.

Now define a local $2$-form
$\bar{\fund}_U\ug s_U^*i^*\fund$ on $U$. Since vertical vector fields are \fund-orthogonal
to any vector field on
$\pi^{-1}(U)$, this definition implies that  for any pair
$(\bar{V},\bar{W})$ of vector fields on $U$
\[
\begin{split}
\bar{\fund}_U(\bar{V},\bar{W})&=s_U^*i^*\fund(\bar{V},\bar{W})\\
&=i^*\fund(V+\nu_V,W+\nu_W)\\
&=i^*\fund(V,W).
\end{split}
\]

Since $i^*\fund$ is compatible with $J$ and
positive, the local form $\bar{\fund}_U$
easily turns out to be
compatible with $\bar{J}$, since
\[
\begin{split}
\bar{\fund}_U(\bar{J}(\bar{V}),\bar{J}(\bar{W}))&=s_U^*i^*\fund(\bar{J}(\bar{V}),\bar{J}(\bar{W}))\\
&=i^*\fund(ds_U(\pi_*{J}({V})),ds_U(\pi_*{J}({W})))\\
&=i^*\fund(J(V)+\nu_{J(V)},J(W)+\nu_{J(W)})\\
&=i^*\fund(J(V),J(W))\\
&=i^*\fund(V,W)\\
&=\bar{\fund}_U(\bar{V},\bar{W})
\end{split}
\]
and the same way one shows that $\bar{\fund}_U$ is positive.

Denote by
$\bar{g}_U$ the corresponding local Hermitian metric, which is
then locally conformal K\"ahler.

We want now to show that $\pi^*\bar{\fund}_U$ is conformally
equivalent to $i^*\fund$ on $\pi^{-1}(U)$.

So consider a pair of generic (that is, non necessarily projectable) vector fields
$(\tilde{V},\tilde{W})$ on $\pi^{-1}(U)$.  For any
$x\in\pi^{-1}(U)$ denote by
${V}^x$ the projectable vector field such that $V^x_x$ coincide with
$\tilde{V}_x$, that is $V^x_y\ug(h_{x,y}^{-1})_*\tilde{V}_{h_x^{-1}s_U(\pi(y))}$,
where by $h_{x,y}$ we denote the element of $G$ that takes
$y$ in $h_x^{-1}s_U(\pi(y))$.
Similarly define ${W}^x$, and call
$(\bar{V}^x,\bar{W}^x)$ the projected vector fields on
$U$.
We then have
\[
\begin{split}
\pi^*\bar{\fund}_U(\tilde{V}_x,\tilde{W}_x)&=\bar{\fund}_U(\pi_*\tilde{V}^x_{x},\pi_*\tilde{W}^x_{x})\\
&=\bar{\fund}_U(\bar{V}^x_{\pi(x)},\bar{W}^x_{\pi(x)})\\
&=i^*\fund(V^x_{s_U(\pi(x))},W^x_{s_U(\pi(x))}).
\end{split}
\]
By evaluating  the projectable vector field
$V^x$ in the point $y=s_U(\pi(x))$ one obtains the following
$$
\begin{array}{ll}
\pi^*\bar{\fund}_U(\tilde{V}_x,\tilde{W}_x)&=i^*\fund((h_x)_*\tilde{V}_x,(h_x)_*\tilde{W}_x)\\
&=h_x^*i^*\fund(\tilde{V}_x,\tilde{W}_x).\\
\end{array}
$$
Now remark that $h_x$ is a conformal map, hence there exists a
smooth function $\alpha_x$ such that
$h_x^*i^*\fund(\tilde{V}_x,\tilde{W}_x)=
\alpha_x(x)i^*\fund(\tilde{V}_x,\tilde{W}_x)$. But by construction
the function $x\mapsto\alpha_x(x)$ is smooth, so
the two 2-forms are conformally equivalent.

Then, if $U, U'\in\mathcal{U}$ overlap, we obtain
on their intersection
that $\bar{\fund}_U$ is conformally
equivalent to $\bar{\fund}_{U'}$:
\[
\bar{\fund}_{U'}=s_{U'}^*i^*\fund\sim
s_{U'}^*\pi^*\bar{\fund}_U=\bar{\fund}_U.
\]
We use a partition of
unity $\{\rho_U\}$ to glue all together these local forms,\label{partitionpage}
obtaining a global $2$-form
$$
\bar{\fund}=\sum_{U\in\mathcal{U}}\rho_U\bar{\fund}_U
$$
on $\moment^{-1}(0)/G$ which, by Lemma~\ref{incollare}, is locally conformal
to any $\bar{\fund}_U$.

 This implies that $\bar{\fund}$ is
still compatible with $\bar{J}$ and positive, and
therefore induces a global Hermitian metric
$\bar{g}$ on $\moment^{-1}(0)/G$ which is locally conformal
K\"ahler because it is
locally conformal to the locally conformal
K\"ahler metrics $\bar{g}_U$ on $U$. This ends the existence
part.

If $g'$ is any locally conformal K\"ahler metric on
$\moment^{-1}(0)/G$ such that $\pi^*g'\sim i^*g$, then for any $x\in\moment^{-1}(0)/G$
on $U_x\subset U$ we obtain
$g'|_{U_x}=s_U^*\pi^*g'|_{U_x}\sim s_U^*i^*g|_{U_x}=\bar{g}_U|_{U_x}\sim
\bar{g}|_{U_x}$. So the globally defined metrics $g$ and $g'$,
being locally conformal, are in fact conformal, by Lemma~\ref{incollare}.
The claim then follows.
\end{D}

\begin{oss}
If $\moment^{-1}(0)/G$ has real dimension two then reduction
equips it with a complex structure and a conformal family of K\"ahler
metrics.
\end{oss}

\begin{oss}\label{crremark}
Let us note by passing that the zero level set offers a natural example of
CR-subman\-ifold of $M$ (see \cite{DrOLCK}). Indeed, the tangent space in
each point splits as a direct
orthogonal sum of a $J$-invariant and a
$J$-anti-invariant distribution:
$T_x(\moment^{-1}(0))=E_x\oplus\liealg{g}(x)$.
 A result of D.~Blair and B.~Y.~Chen
states that the anti-invariant distribution of a CR-submanifold in a
locally conformal K\"ahler manifold is integrable. In our case, this is
trivially true because the anti-invariant distribution is just a copy of the
Lie algebra of $G$.
\end{oss}

\section{Conformal setting and the
Biquard-Gauduchon construction}\label{biqgau}

In defining the reduced locally conformal K\"ahler structure on
$\moment^{-1}(0)/G$ we used a specific metric in the conformal class
$[g]$, to obtain a conformal class $[\bar{g}]$.
In this section we present a more intrinsic construction for the locally
conformal K\"ahler reduction,  due to O.~Biquard and P.~Gauduchon, which makes use of the
language of conformal
geometry.
To this aim we mainly fill in details and reorganize material contained
in \cite{CaPEWG} and in the unpublished paper \cite{BiGCML}.

Moreover we prove that the two constructions are in fact the same, by
showing in Lemma~\ref{lemmalink} and its consequences
the correspondence between representatives and intrinsic
objects.

Let $V$ be a real $n$-dimensional vector space, and $t$ a real
number. The $1$-dimensional vector space $L^t_V$ of
{\em densities of weight $t$} on $V$ is the vector space of maps \map{l}{(\Lambda^n V)\setminus
0}{\R} satisfying $l(\lambda w)=|\lambda|^{-t/n}l(w)$ if
$\lambda\in\R\setminus 0$ and $w\in(\Lambda^n V)\setminus 0$.
We say that a
density $l$ is {\em positive}\/ if it takes only positive real
values. For positive integers $t$ we have $L^t_V=L^1_V\otimes\dots\otimes
L^1_V$ and for negative integers $t$ we have
$L^t_V=(L^1_V)^*\otimes\dots\otimes (L^1_V)^*$. Thus, given an element
$l$ of $L^1_V$, we denote by $l^t$ the corresponding element of $L^t_V$ under
these canonical identifications, for any $t$ integer.

Remark that $\Lambda^{n+d}(V\oplus\R^d)\iso\Lambda^n V$, and this
gives a canonical isomorphism between $L^t_{V\oplus\R^d}$ and $L^t_V$:
\begin{equation}\label{isodensita}
l\in L^t_V \mapsto \text{\rm sgn}(l) l^{\frac{n}{n+d}}\in L^t_{V\oplus\R^d}.
\end{equation}

To any Euclidean metric
$g$ on the vector space $V$ we associate the positive element $l_g^t$ of $L^t_V$
which sends the length-one element of $(\Lambda^nV)\setminus
0$ to $1$. Then under a homothety $e^\alpha g$ of the metric we
have $l_{e^\alpha g}^t=e^{-t\alpha/2}l_g^t$, and the positive definite
element $g\otimes l_g^2$ of
$S^2V\otimes L^2_V$ only depends  on the homothety class $c$ of $g$.

Conversely, given an element $c$ of $S^2V\otimes L^2_V$, we can
associate to any positive element $l$ of $L^1_V$ the
element $c\otimes l^{-2}$ of $S^2V\otimes L^2_V\otimes
L^{-2}_V=S^2V$, and if $c$ is positive definite so is $c\otimes l^{-2}$,
which therefore defines a Euclidean metric on $V$. If, moreover,
$c$ satisfies the normalization condition $l^2_{c\otimes
l^{-2}}=l^2$ for one (and hence for all) positive element
$l$ of $L^1_V$, then the correspondence between such $c$'s and the homothety classes of
$g$ is bijective.

For any vector bundle $E\rightarrow M$, define the associated
{\em density line bundle} $L^t_E\rightarrow M$ as the bundle
whose fiber over $x\in M$ is the $1$-dimensional vector space
$L^t_{E_x}$. If $n$ is the rank of $E$, then $L^t_E$ can be
globally defined as the fibred product
$\pr{E}\times_{G}L^t_{\R^n}$, where $\pr{E}$ denotes the principal bundle
associated to $E$ with structure group $G\subset\gl{n}$, and an element $A$ of $G$
acts on $L^t_{\R^n}$ by multiplication by $|\det A|^{t/n}$.
Remark that, in particular, $L^t_E$ has the same principal bundle
as $E$, for any $t\in\R$.

The above construction identifies conformal classes of metrics on
$E$\/ with normalized positive defined sections of $S^2E\otimes L^2_E$. In
particular, if $E=TM$, the conformal class of a Riemannian metric
can be thought of as a normalized positive defined section $c$ of $S^2M\otimes
L^2_M$, where we denote $L^t_{TM}$ by $L^t_M$.

A trivialization (usually positive) of $L^1_M$ is called a {\em
gauge} or also a {\em length scale}.

This way, on a conformal manifold $(M,c)$, we have a Riemannian metric
whenever we fix a gauge. As a terminology, instead of saying ``\dots take a gauge $l$,
and let $g\ug c\otimes l^{-2}$\dots'' we shall say ``\dots let $g$ be a metric
in the conformal class $c$\dots''.

Since a connection on $M$ means a connection on $\gl{M}$ and
$\gl{M}$ is also the principal bundle of $L^t_M$, a connection on
$M$
 induces a
connection on $L^t_M$, for any $t\in\R$. Vice versa, suppose a
connection $\nabla$ on $L^1_M$ is given. Then we can use a conformal
version of the
six-terms formula to define a connection on
$M$, still denoted by $\nabla$, which is compatible with $c$:
\begin{equation}\label{conformalsix}
2c(\nabla_X Y,Z)=\nabla_X c(Y,Z)+\nabla_Y c(X,Z)-\nabla_Z
c(X,Y)+c([X,Y],Z)-c([X,Z],Y)-c([Y,Z],X),
\end{equation}
where both members are sections of $L^2_M$.

This way one proves
the fundamental
theorem of conformal geometry:
\begin{teo}[Weyl]
Let $(M,c)$ be a conformal manifold.
There is an affine bijection between connections on $L^1_M$ and
torsion-free connections on $M$ preserving $c$.
\end{teo}

Torsion-free compatible connections on a conformal manifold are called Weyl connections.
In contrast with the Riemannian case, the previous theorem says in particular
that on a conformal
 manifold there is not a uniquely defined torsion-free
compatible connection.

In this setting a {\em conformal almost-Hermitian manifold} is a conformal
manifold $(M,c)$ together with an almost-complex structure $J$ on $M$
compatible with one (and hence with all) metric in the conformal class.

Let $(M,c,J)$ be a conformal almost-Hermitian manifold.
We then have a non-degenerate
fundamental form $\fund$ taking values in $L^2_M$, that is, $\fund(X,Y)\ug
c(JX,Y)\in\Gamma(L^2_M)$, for $X,Y\in\Chi(M)$.
For any metric
$g$ defining $c$, with corresponding fundamental form $\fund_g$, we have
$\fund=\fund_g\otimes l^2_g$.
The notion of Lee form $\lee_g$ of the almost-Hermitian metric $g$ on $(M,J)$ is clearly dependent on
the metric,
but a straightforward computation shows that the connection $\nabla$
on $L^1_M$ given by $\nabla_Xl_g\ug(-1/2)\lee_g(X)l_g$
does not depend on the choice of $g$ in the conformal class $c$.

The fundamental theorem of conformal geometry
gives then a torsion-free compatible connection on $M$, which is called the {\em canonical Weyl
connection} of the conformal almost-Hermitian manifold $(M,c,J)$. We
denote simply by $\nabla$ this connection on $M$, and we use the
same symbol for the induced connection on $L^t_M$, for any
$t\in\R$. In particular, the constant $-1/2$ in the definition of $\nabla$
was chosen in order that $\nabla l_g^2=-\lee_g\otimes
l_g^2$.

Thus, given any $L^t_M$-valued tensor $\psi$ on a conformal almost-Hermitian manifold,
we can differentiate it with respect to the canonical Weyl
connection, and any choice of a metric $g$ in the conformal class $c$
gives a corresponding real valued tensor $\psi_g$. The following
Lemma links this intrinsic point of view with the gauge-dependant
setting of almost-Hermitian manifolds. We state it only for
$L^2_M$-valued differential forms, because this is the only case
we need.

\begin{lem}[Equivalence lemma]\label{lemmalink}
Let $(M,c,J)$ be a conformal almost-Hermitian manifold, with canonical Weyl connection
$\nabla$. Let $\psi$ be a $p$-form taking values in $L^2_M$.
Then for any metric $g$ in the conformal class $c$ we have
\[
d^\nabla\psi=d^{\lee_g}\psi_g\otimes l_g^2.
\]
\end{lem}

\begin{D}
\begin{equation*}
\begin{split}
d^\nabla\psi&=d^\nabla(\psi_g\otimes l^2_g)=d\psi_g\otimes
l^2_g+(-1)^{|\psi_g|}\psi_g\wedge\nabla l_g^2\\
&=d\psi_g\otimes
l_g^2-(-1)^{|\psi_g|}\psi_g\wedge\lee_g \otimes l_g^2=d\psi_g\otimes
l_g^2-\lee_g\wedge\psi_g\otimes l_g^2=d^{\lee_g}\psi_g\otimes
l_g^2.
\end{split}
\end{equation*}
\end{D}

Using the equivalence Lemma we obtain in particular
\begin{equation}\label{denablaomega}
d^\nabla\fund=d^{\lee_g}\fund_g\otimes
l_g^2.
\end{equation}

Since the Weyl connection is compatible with $c$, we have also
\begin{equation}\label{nablac}
0=\nabla c=\nabla(g\otimes l_g^2)=\nabla g\otimes l_g^2+g\otimes
\nabla l_g^2=\nabla g\otimes l_g^2-g\otimes\lee_g\otimes
l_g^2=(\nabla g-\lee_g\otimes g)\otimes
l_g^2.
\end{equation}

\begin{teo}\label{almost}
Let $(M,c,J)$ be a conformal almost-Hermitian
manifold, and let  $\nabla$  be the canonical Weyl connection. Let
$g$ be any metric in the conformal class $c$.
Then:
\begin{enumerate}
\item $\nabla$ preserves $J$ if and only if $J$ is integrable
and $(d\fund_g)_0=0$;
\item\label{flat} the curvature $R^{\nabla}=\nabla_{[X,Y]}-[\nabla_X,\nabla_Y]$ of\, $\nabla$ is given
by $R^{\nabla}l_g^2=d\lee_g\otimes l_g^2$.
\end{enumerate}
\end{teo}

\begin{D}
For any complex connection $\nabla$ the following formula holds, linking the torsion $T$
of $\nabla$ with the torsion $N$ of $J$:
\[
T(JX,JY)-J(T(JX,Y))-J(T(X,JY))-T(X,Y)=-N(X,Y).
\]
Since Weyl connections are torsion free, if we find any complex Weyl
connection then  $J$ is
integrable. We want to show that, if the canonical Weyl connection
is complex, then also $(d\fund_g)_0=0$.
Denote by
$A$ the alternation operator and by $C$ the contraction such that $\fund_g=C(J\otimes g)$, then
\begin{equation*}
\begin{split}
d\fund_g&=A(\nabla\fund_g)=A(\nabla C(J\otimes g))=A(C(J\otimes\nabla
g))=A(C(J\otimes\lee_g\otimes g))\\
&=A(\lee_g\otimes C(J\otimes
g))=A(\lee_g\otimes\fund_g)=\lee_g\wedge\fund_g,
\end{split}
\end{equation*}
where we have used formula \eqref{nablac} to obtain $\nabla g=\lee_g\otimes
g$.

Suppose now that $(d\fund_g)_0=0$ and that $J$ is integrable. Then
using the conformal six-terms formula \eqref{conformalsix}
we obtain the following
conformal version of a classical formula in Hermitian geometry (see \cite[p.~148]{KoNFD2}):
\[
4 c((\nabla_X J)Y,Z)=6 d^\nabla \fund(X,JY,JZ)-6 d^\nabla
\fund(X,Y,Z),
\]
and this shows that $c((\nabla_X J)Y,Z)=0$ if $d^\nabla\fund=0$.
But this last condition is equivalent, by formula
\eqref{denablaomega}, to  $(d\fund_g)_0=0$, and the
claim then follows from the non-degeneracy
of $c$. As for the curvature $R^\nabla$ of $\nabla$, using equivalence Lemma we obtain
\[
R^\nabla l_g^2=-d^\nabla(\nabla l_g^2)=-d^\nabla(-\lee_g\otimes
l_g^2)=d^{\lee_g}\lee_g\otimes l_g^2=d\lee_g\otimes
l_g^2.
\]
\end{D}

Since a locally conformal K\"ahler manifold is a conformal
Hermitian manifold $(M,c,J)$ such that
$(d\fund_g)_0=0$ and
$d\lee_g=0$, for one (and then for all)
choice of metric $g$ in the conformal class $c$ (compare with formula
\eqref{itb}),
we can give the following intrinsic
characterization of locally conformal K\"ahler
manifolds:

\begin{cor}\label{caratterizzazione}
Let $(M,c,J)$ be a conformal Hermitian manifold. Denote by $\fund$
the $L^2_M$-valued fundamental form, and
let $\nabla$ be the canonical
Weyl connection. Then $(M,c,J)$ is
locally conformal K\"ahler if and only if
$\nabla$ is flat and $\fund$ is $d^\nabla$-closed.
\end{cor}

%\begin{D}
%We have already remarked that $d^\nabla \fund=0$ is equivalent to
%$(d\fund_g)_0=0$, for any $g$ in the conformal class $c$.
%Let's compute the curvature $R^\nabla$ of $\nabla$:
%\[
%R^\nabla l_g^2=-d^\nabla(\nabla l_g^2)=-d^\nabla(-\lee_g\otimes
%l_g^2)=d\lee_g\otimes l_g^2-\lee_g\wedge\lee_g\otimes l_g^2=d\lee_g\otimes
%l_g^2.
%\]
%Then $R^\nabla=0$ if and only if the Lee form $\lee_g$ is closed for one (and then for all)
%choice of $g$ in the conformal class $c$, and the claim follows.
%\end{D}

Moreover, Theorem \ref{almost} gives also the following
\begin{cor}\label{weylpreservesJ}
On a locally conformal K\"ahler manifold the canonical Weyl
connection preserves the complex structure.
\end{cor}

Unless otherwise stated, from now on we consider locally conformal K\"ahler
manifolds $(M,c,J)$.

A locally conformal K\"ahler manifold $(M,c,J)$ comes then
naturally equipped with a closed $2$-form $\fund$, the only
difference from the K\"ahler case being that $\fund$ now takes
values in $L^2_M$. We want go further with this analogy.

Define the pairing
%\map{\flat}{\Chi(M)}{\Omega^1(L^2_M)} by $\flat X\ug \iota_X\fund$,
%and the corresponding dual
\map{\sharp}{\Omega^1(L^2_M)}{\Chi(M)} by
$\iota_{\sharp\alpha}\fund=\alpha$, and use it to
define a Poisson bracket on
$\Gamma(L^2_M)$ by
%\begin{equation*}
$\{f_1,f_2\}\ug\fund(\sharp \nabla f_1,\sharp \nabla f_2)$.
%\end{equation*}
Using Lemma \ref{lemmalink} and  formula
\eqref{lieid}, one shows the
relation
\begin{equation}\label{lieidconf}
\{\{f_1,f_2\},f_3\}+\{\{f_2,f_3\},f_1\}+\{\{f_3,f_1\},f_2\}
=d^\nabla\fund(\sharp \nabla f_1,\sharp \nabla f_2,\sharp \nabla
f_3)=0,
\end{equation}
proving that this bracket turns $\Gamma(L^2_M)$ into a Lie
algebra. Remark that, as formula \eqref{lieid}, the first equality
in \eqref{lieidconf}  holds generally on conformal
almost-Hermitian manifolds such that the canonical Weyl connection is
flat.

We finally describe the intrinsic version of $\aut{M}$. If
$l$ is a section of $L^t_M$ and $h$ is a diffeomorphism of $M$, then the section
$h^*l$ of $L^t_M$ is given
by
\begin{equation}\label{azionesuL}
(h^*l)_x\ug l_{h(x)}\circ (h_*)_x,
\end{equation}
 that is, if
$x\in M$ and $w\in \Lambda^n(T_xM)\setminus0$, we have
$(h^*l)_x(w)\ug l_{h(x)}((h_{*})_xw)$. Recall that
for any $x$ the differential induces the map
$\map{(h_*)_x}{\Lambda^n(T_xM)}{\Lambda^n(T_{h(x)}M)}$ which is in
fact a linear map between $1$-dimensional vector spaces. Whenever a
metric $g$ is fixed, a trivialization $w^g$
of $\Lambda^n(TM)$ associating to $x$ the length-one element
$w_x^g$ is defined, hence one can associate to any diffeomorphism
$h$ a never-vanishing smooth function $d_h^g$ defined by
\[
h_*w^g_x=d_h^g(x)w^g_{h(x)},
\]
so the following derivation rule holds for $l_g^t$:
\[
\begin{split}
(h^*l_g^t)_x(w^g_x)&=(l_g^t)_{h(x)}(h_*w^g_x)\\
&=(l_g^t)_{h(x)}(d^g_h(x) w^g_{h(x)})\\
&=|d_h^g(x)|^{-\frac{t}{n}}(l_g^t)_{h(x)}(w^g_{h(x)})=|d_h^g(x)|^{-\frac{t}{n}},
\end{split}
\]
that is, in short, $h^*l_g^t=|d_h^g|^{-\frac{t}{n}}l_g^t$.
%
%=|\det (h_{*})_x|^{-\frac{1}{n}}l_{h(x)}(w)$, that is, in short,
%$h^*l=|\det h_*|^{-\frac{1}{n}}l$.

For any diffeomorphism of $M$ we
then define $h^*c$ in the obvious way, that is,
$h^*c=h^*g\otimes h^*l^2_g$. Since
\[
d_h^{e^fg}=e^{(n/2)f\circ h}e^{(-n/2)f}d_h^g,
\]
this definition does not depend on
the choice of the gauge $g$, and
gives the intrinsic notion of $\aut{M}$ as follows.
\begin{pro}
A diffeomorphism $h$ of a conformal manifold $(M,c)$ preserves $c$
if and only if it is a conformal transformation of one (and hence
of all) metric $g$ in the conformal class $c$.
\end{pro}

\begin{D}
Indeed,  $h^*g=e^\alpha g$ implies $d^g_h=e^{n\alpha/2}$, so
$h^*l_g^2=e^{-\alpha}l_g^2$, and then $h^*c=h^*(g\otimes l_g^2)=h^*g\otimes h^*l_g^2=g\otimes
l_g^2= c$. Vice versa $h^*c=c$ implies $h^*g\otimes h^*l_g^2=g\otimes
l_g^2$, hence $(|d^g_h|^{-\frac{2}{n}}h^*g)\otimes
l_g^2=g\otimes l_g^2$, that is $h^*g=|d^g_h|^{\frac{2}{n}}g$.
\end{D}

\begin{lem}\label{invarianza}
The Weyl connection of a conformal almost-Hermitian manifold $(M,c,J)$ is
invariant for \aut{M}, that is, $h_*\nabla_VW=\nabla_VW$ whenever
$h_*V=V$ and $h_*W=W$.
\end{lem}

\begin{D}
This is because \aut{M} preserves $c$ and $J$, and
$\nabla$ is defined just using these ingredients. More formally,
we want to show that, if $h\in\aut{M}$ and $V$, $W$, $Z$ are
$h$-invariant vector fields, then
$c(\nabla_VW,Z)=c(h_*\nabla_VW,Z)$. But we have
\[
c(\nabla_VW,Z)=(h^*c)(\nabla_VW,Z)=h^*(c(h_*\nabla_VW,h_*Z))=h^*(c(h_*\nabla_VW,Z)),
\]
where we used the general property that if $\psi$ is any tensor field
of type $(r,0)$ and  $X_1,\dots,X_r$ are vector fields, then
$h^*(\psi(h_*X_1,\dots,h_*X_r))=(h^*\psi)(X_1,\dots,X_r)$.
We are therefore only left to show that $c(\nabla_VW,Z)$ is
$h$-invariant  for all $h\in\aut{M}$, that is, we are left to show that
the second side of the conformal six-terms formula \eqref{conformalsix}
is $h$-invariant for all $h\in\aut{M}$. But it turns out that each
summand of  \eqref{conformalsix} is $h$-invariant.
We show this only on its first and fourth summand,
the others being similar: the first summand
\begin{equation*}
\begin{split}
h^*\nabla_Vc(W,Z)&=h^*\nabla_V(g(W,Z)l^2_g)\\
&=h^*(Vg(W,Z)l^2_g)+h^*(g(W,Z)\nabla_Vl^2_g)\\
&=h^*Vg(W,Z)h^*l^2_g-h^*(g(W,Z))h^*(\lee_g(V))h^*l^2_g\\
&=V((h^*g)(W,Z))l^2_{h^*g}-(h^*g)(W,Z)\lee_{h^*g}(V)l^2_{h^*g}\\
&=\nabla_V((h^*g)(W,Z)l^2_{h^*g})=\nabla_Vc(W,Z),
\end{split}
\end{equation*}
where we have used that $V$ and $h$ commute on $C^\infty(M)$,
since $V$ is $h$-invariant, that $h^*l^2_g=l^2_{h^*g}$ and that
$h^*\lee_g=\lee_{h^*g}$. The fourth summand is
\begin{equation*}
\begin{split}
h^*(c([V,W],Z))=h^*(g([V,W],Z)l^2_g)=h^*(g([V,W],Z))h^*l^2_g=(h^*g)([V,W],Z)l^2_{h^*g}=c([V,W],Z),
\end{split}
\end{equation*}
where we have used the already cited properties and that the Lie
bracket of invariant vector fields is invariant.
\end{D}

\begin{cor}\label{corinvarianza}
Let $(M,c,J)$ be a conformal almost-Hermitian manifold with
Weyl connection $\nabla$, and let
$G\subset\aut{M}$. If $V$, $W$, $Z$ are $G$-invariant vector fields
on $M$, then $c(\nabla_VW,Z)$ is $G$-invariant.
\end{cor}

Let $G$ be a Lie subgroup of $\aut{M}$, as in section~\ref{momentmap}.
The momentum map can then be defined as a homomorphism of Lie
algebras
\map{\momenthat^\cdot}{\liealg{g}}{\Gamma(L^2_M)}
such that
$\iota_X\fund=d^\nabla\momenthat^X$. We also denote by $\moment$ the
corresponding element of $\Gamma(\liealg{g}^*\otimes L^2_M)$ given by
$\langle\moment(x),X\rangle=\momenthat^X(x)$, carets denoting the evaluation.

\begin{oss}
In \cite{BiGCML}
the existence of such a homomorphism of Lie algebras is shown to imply
the condition
\[
-\frac{1}{2}\lee_g(X)+\frac{1}{n}\mathop{\mathrm{div}_g} X=0
\]
on any fundamental vector field $X$. This is
equivalent to the condition
\[
\lie_X\fund_g-\lee_g(X)\fund_g=0
\]
one finds in~\cite{HaRRLC}, since $\lie_X\fund_g=((2/n)\mathop{\mathrm{div}_g}
X)\fund_g$.
\end{oss}

If we choose a metric $g$ in the conformal class $c$, then $\momenthat^X=\momenthat_g^X
l_g^2$, where \map{\momenthat_g^\cdot}{\liealg{g}}{C^\infty(M)}.

\begin{teo}\label{teomoment}
The map \map{\momenthat^\cdot}{\liealg{g}}{\Gamma(L^2_M)} is a
momentum map if and only if \map{\momenthat_g^\cdot}{\liealg{g}}{C^\infty(M)} is
a momentum map as in section \ref{momentmap}.
\end{teo}

\begin{D}
Use Lemma \ref{lemmalink} to compute $d^\nabla\momenthat^X$ with respect to the fixed gauge:
\[
d^\nabla\momenthat^X=d^{\lee_g}\momenthat_g^X\otimes l_g^2,
\]
so that $d^\nabla\momenthat^X=\iota_X(\fund_g\otimes l_g^2)=\iota_X\fund_g\otimes
l_g^2$ if and only if $\iota_X\fund_g=d^{\lee_g}\momenthat_g^X$.
We then have to check that $\momenthat$ is a Lie algebra
homomorphism if and only if $\momenthat_g$ is. But this is a
direct consequence of Lemma \ref{lemmalink}, and of the fact that
$\sharp\alpha=\sharp_g\alpha_g$:
\[
\{f_1,f_2\}=\fund(\sharp \nabla f_1,\sharp \nabla f_2)=\fund_g(\sharp \nabla f_1,\sharp \nabla
f_2)\otimes l_g^2=\fund_g(\sharp_g d^{\lee_g} f_{1,g},\sharp_g d^{\lee_g} f_{2,g})\otimes
l_g^2=\{f_{1,g},f_{2,g}\}_g\otimes l^2_g.
\]
\end{D}

\begin{oss}
The previous theorem allows using all proofs of section
\ref{momentmap} as proofs in this conformal setting, just fixing a
gauge. In particular, the zero set $\moment^{-1}(0)$, where $0$
denotes the zero section of $\liealg{g}^*\otimes L^2_M$, is the zero set
of any $\moment_g$, and it is therefore closed with respect to the
action of $G$ and coisotropic with respect to $\fund$. Moreover,
the assumption of $0$ being a regular value for $\moment_g$
translates into the assumption that the zero section be transverse
to $\moment$, and under this assumption the
isotropic foliation is given exactly by fundamental vector fields
$\liealg{g}$.
\end{oss}

\begin{teo}[Biquard \& Gauduchon, \cite{BiGCML}]
Let $(M,c,J)$ be a locally conformal K\"ahler manifold.
Let $G$ be a Lie subgroup of\/ $\aut{M}$ whose action admits a momentum
map \map{\momenthat}{\liealg{g}}{\Gamma(L^2_M)}. Suppose that $G$ acts freely and
properly on $\moment^{-1}(0)$, $0$ denoting the zero section of\, $\liealg{g}^*\otimes
L^2_M$, and suppose that $\moment$ is transverse to this zero
section.
Then
there exists a locally conformal K\"ahler structure $(\bar{c},\bar{J})$ on
$\moment^{-1}(0)/G$.
% uniquely determined by the
%condition $\pi^* \bar{c}=i^*c$,
%where $i$ denotes the inclusion of $\moment^{-1}(0)$ into $M$ and
%$\pi$ denotes the projection of $\moment^{-1}(0)$ onto its
%quotient.
\end{teo}

\begin{D}
Due to Lemma \ref{lemmalink} and to Theorem \ref{teomoment}, this
theorem can be viewed at as a translation of Theorem \ref{teonostro}
in the conformal language. From this point of view, the theorem
was already proved.

We want here to give an intrinsic proof, using the
characterization of locally conformal K\"ahler manifolds given by
corollary~\ref{caratterizzazione}.

Take the $c$-orthogonal decomposition $T_xM=E_x\oplus\liealg{g}(x)\oplus
J\liealg{g}(x)$, where $E_x$ is the $c$-orthogonal complement of $\liealg{g}(x)$ in
$T_x(\moment^{-1}(0))$. We obtain a vector bundle $E\rightarrow
\moment^{-1}(0)$ of rank $n-2\dim G$.

First we need to relate $L^t_{\moment^{-1}(0)/G}$ with $L^t_E$.
Remark that ${E/G}\rightarrow{\moment^{-1}(0)/G}$ is isomorphic as a bundle to the tangent bundle
of
${\moment^{-1}(0)/G}$, by means of $\pi_*|_E$. On
its side
$L^t_E/G$ is isomorphic to $GL(E/G)\times_{GL(n-2\dim G)} L^t_{\R^{n-2\dim G}}$, since the actions of $G$ and of
$GL(n-2\dim G)$ on $GL(E)$ commute, that is, if $g\in G$,
$\gamma\in GL(n-2\dim G)$ and $p\in GL(E)$, then
$g_*(p\gamma)=(g_*p)\gamma$.
This means that $L^t_{\moment^{-1}(0)/G}$ is isomorphic to
$L^t_E/G$, the isomorphism being explicitly given by sending
an element
$l$ of $L^2_{\moment^{-1}(0)/G,\bar x}$ to
$[l\circ \pi_{*,x}]$,  where $\pi(x)=\bar x$, and the action of $G$
on $L^t_E$ being given  by \eqref{azionesuL}.

Now remark that the canonical splitting
$TM=E\oplus\liealg{g}\oplus
J\liealg{g}$ gives an isomorphism of $L^2_M|_{\moment^{-1}(0)}$
with $L^2_E$, by formula \eqref{isodensita}, and this isomorphism
is $G$-equivariant.

We therefore think of elements of $L^2_{\moment^{-1}(0)/G}$ as
equivalence classes of elements of $L^2_M|_{\moment^{-1}(0)}$.

During the proof of this theorem, we denote by $\bar{V},\bar{W},\dots$
vector fields on ${\moment^{-1}(0)/G}$, and by $V,W,\dots$ their lifts to
$E$. Note that $V,W,\dots$ are $G$-invariant vector fields.

Define
$\bar{c}(\bar{V},\bar{W})$ to be the projection to $L^2_M|_{\moment^{-1}(0)}/G$ of
the section $c(V,W)$, that is
\[
(\bar{c}(\bar{V},\bar{W}))_{\bar{x}}\ug [c(V,W)_{x}]\in
(L_M^2|_{\moment^{-1}(0)})_x/G\simeq (L_{\moment^{-1}(0)/G}^2)_{\bar{x}}
\]
where $x$ is an element in $\pi^{-1}(\bar{x})$. The choice of $x$
is irrelevant, since
$h^*(c(V,W))=h^*(c(h_*V,h_*W))=(h^*c)(V,W)=c(V,W)$.

We have thus defined an almost-Hermitian conformal manifold
$(\moment^{-1}(0)/G,\bar{c},\bar{J})$. In order to show that it is
locally conformal K\"ahler we compute its canonical Weyl
connection, and then use corollary \ref{caratterizzazione}.
%introduce a connection
%$\bar{\nabla}$ that we will prove to be the Weyl connection of
%$(\moment^{-1}(0)/G,\bar{c},\bar{J})$.

Let
$\nabla^E$ be the orthogonal projection of $\nabla$ from $T(\moment^{-1}(0))$ to $E$.
Since by Lemma \ref{invarianza} the
Weyl connection $\nabla$ is invariant for
$\aut{M}$, we have  that $\nabla^E_VW$ is a projectable vector
field. Define
\begin{equation}\label{barnabla}
\bar{\nabla}_{\bar{V}}\bar{W}\ug \pi_* \nabla^E_VW.
\end{equation}
The torsion $T^{\bar\nabla}_{\bar V,\bar W}$ of $\bar\nabla$ is just
$\pi_*T^{\nabla^E}_{V,W}=0$. Moreover,
$\bar\nabla$ is compatible with $\bar J$:
\[
\begin{split}
(\bar\nabla_{\bar V}\bar J)\bar W
&=\bar\nabla_{\bar V}(\bar J\bar W)-\bar J\bar\nabla_{\bar V}\bar
W\\
&=\pi_*\nabla^E_{V}(JW)-\bar J\pi_*\nabla^E_{V} W\\
&=\pi_*(\nabla^E_{V}(JW)-J\nabla^E_{V} W)
=\pi_*(\nabla_{V}J)^EW=0.
\end{split}
\]
Eventually, Theorem \ref{almost} proves that $\bar J$
is integrable.

Look at the Weyl connection $\nabla$ on $L^2_M$ as a map
\map{\nabla_V}{\Gamma(L^2_M|_{\moment^{-1}(0)})}{\Gamma(L^2_M|_{\moment^{-1}(0)})},
and remark that the \aut{M}-invariance of $V$ implies that
$\nabla_V$ is $G$-equivariant, thus defines a connection on
$L^2_{\moment^{-1}(0)/G}$. We denote it again by $\bar\nabla$:
\[
\bar\nabla_{\bar V}[l]\ug[\nabla_Vl]\in
L^2_M|_{\moment^{-1}(0)}/G.
\]
%}

Using the conformal six-terms formula
\eqref{conformalsix} and corollary \ref{corinvarianza}, we see
that the connection $\bar\nabla$ on $\moment^{-1}(0)/G$ defined by
\eqref{barnabla} is the associated Weyl connection, which is therefore
the canonical Weyl connection of $(\moment^{-1}(0)/G,\bar c,\bar
J)$.

The curvature $R^{\bar\nabla}$ is given by
\[
R^{\bar\nabla}_{\bar V,\bar W}[l]=-d^{\bar \nabla}\bar\nabla[l](\bar V,\bar W)=-[d^\nabla\nabla
l(V,W)]=[R^\nabla_{V,W}]=0.
\]

Finally, denoting by $\bar\fund$ the
$L^2_{\moment^{-1}(0)/G}$-valued fundamental form of
$(\moment^{-1}(0)/G,\bar c,\bar J)$, we have
$\pi^*(d^{\bar\nabla}\bar\fund)=d^\nabla\fund=0$, thus $d^{\bar\nabla}\bar\fund=0$,
and corollary \ref{caratterizzazione} says that $(\moment^{-1}(0)/G,\bar c,\bar J)$ is
locally conformal K\"ahler.
\end{D}

\section{Compatibility with K\"ahler reduction}\label{trivialclass}

In this section we analyze the relation between locally conformal
K\"ahler reduction of a manifold and K\"ahler reduction of a
covering. We refer to \cite{FutKEM} for the K\"ahler reduction.

As a  first step we show that
the two notions of reduction on globally conformal K\"ahler manifolds are compatible.

\begin{pro}\label{globally} Let $(M,[g],J)$ be a globally conformal K\"ahler
manifold and denote by $g$ a K\"ahler metric. Let $G\subset \aut{M}$ a subgroup satisfying the
hypothesis of the reduction theorem and which moreover is
composed by isometries with respect to $g$. Denote by
$(\moment^{-1}(0)/G,[\bar{g}],\bar{J})$ the reduced locally
conformal K\"ahler manifold.
Then the action of $G$ is Hamiltonian for $g$,
the submanifold $\moment^{-1}(0)$ is the same as in the K\"ahler
reduction
and the conformal class of the reduced K\"ahler metric is $[\bar{g}]$.
So, in particular, the
reduced manifold
is globally conformal K\"ahler.
\end{pro}

\begin{D}
As  the
action of $G$ is twisted Hamiltonian for $[g]$
 Remarks~\ref{conformalhamiltonian} and~\ref{conformalpoisson} imply that
it is Hamiltonian for $g$.
Moreover,
the subspace $\moment^{-1}(0)$ is the same for both
notions.
The construction of the almost-complex structure on the quotient
is the same in the two cases, so $\bar{J}$ is defined.
Denote by $\tilde{\fund}$ the K\"ahler form that the K\"ahler reduction provides on
$\moment^{-1}(0)/G$. Then $\pi^*\tilde{\fund}=
i^*\fund$, so the claim follows by the uniqueness part of the reduction theorem.
\end{D}

\begin{exa}
If $(M,[g],J)$ is a globally conformal K\"ahler
manifold the reduced structure is not necessarily globally conformal
K\"ahler. Actually, any locally conformal K\"ahler manifold
$(M,[g],J)$ can be
seen as a reduction of a globally conformal manifold.
Indeed, consider the universal covering $\tM$ of $M$
equipped with its pulled-back locally conformal K\"ahler
structure, which is globally conformal K\"ahler since \tM\ is
simply connected.
 This covering manifold can be considered to be acted on
by the
discrete group of holomorphic conformal maps $G\ug\pi_1(M)$,
which, having trivial associated infinitesimal
action, is clearly Hamiltonian, with trivial momentum map: hence $\moment^{-1}(0)=\tM$
and $\moment^{-1}(0)/G=M$.
\end{exa}

We now concentrate  our attention to the structure of the universal cover
$\tilde{M}$ of a locally conformal K\"ahler manifold $(M,[g],J)$.

\begin{oss}\label{coveringhomothetic}
The pull-back by the covering map $p$
of any metric  of $[g]$ is globally conformal K\"ahler
since \tM\ is simply connected.
It is easy to show that  on any
complex manifold $Z$ such that $\dim_\C(Z)\geq2$
if two
K\"ahler metrics are conformal then their conformal factor is constant.
 In our case remark that the pull-back of any metric in $[g]$ is conformal to a  K\"ahler
metric $\tilde{g}$ by
\[
\tilde{g}=e^{-\tau}p^*g
\]
where $\tau$ satisfies $d\tau=\lee_{\tilde{g}}=p^*\lee_g$ and is then only
defined up to adding a constant. What is remarkable is that the action
of $\pi_1(M)$ on \tM\ is by homotheties of the K\"ahler metrics  (we fix points in $M$ and in $\tM$ in order
to have this action well-defined).
Moreover any element of \aut{M} lifts to a
homothety of the K\"ahler metrics of
\tM, if $\dim_\C(M)\geq2$. This is in fact an
equivalent definition of locally conformal K\"ahler manifolds
(see~\cite{VaiGHM} and \cite{DrOLCK}).

%We underline this set of facts by saying that
%$\tilde{M}$
%carries a well-defined structure of {\em homothetic
%K\"ahler manifold}.
\end{oss}

With this model in mind, we define a {\em homothetic
K\"ahler manifold} as a triple
$(K,\langle g\rangle,J)$, where $(K,g,J)$ is a K\"ahler manifold
and $\langle g\rangle$ denotes the set of metrics differing from $g$
by multiplication for a positive factor.
We define $\mathcal{H}(K)$ to be the group of
biholomorphisms of $K$ such that $f^*g=\lambda g$, $\lambda\in\R_+$,
and we call such a map a {\em homothety}
 of
$K$ of {\em dilation factor} $\lambda$. The dilation factor does
not depend on the choice of $g$ in $\langle g\rangle$,
so
a homomorphism $\rho$ is defined from $\mathcal{H}(K)$ to $\R_+$
associating to any homothety its dilation factor (see also
\cite{KaOGFC}). Note that
$\ker\rho$ is the subgroup of $\mathcal{H}(K)$ containing the maps
that are isometries of one and then all of the metrics in
$\langle g\rangle$.
If
$K$ is given as a globally conformal K\"ahler manifold $(K,[g],J)$, then $\mathcal{H}(K)$
can be considered as the well-defined subgroup of $\aut{K}$ of
homotheties with respect to the K\"ahler metrics in $[g]$.
We now give a condition for a locally conformal K\"ahler manifold
covered by a globally conformal one to be globally conformal K\"ahler.

\begin{pro}\label{indurre}
Given  a globally conformal K\"ahler manifold $(\tM,[\tilde{g}],J)$
and a subgroup
$\Gamma$
of $\aut{\tM}$ acting freely and properly discontinuously,
the quotient $M\ug \tM/\Gamma$ (with its naturally induced complex
structure)
comes equipped with
a locally conformal K\"ahler structure $[g]$  uniquely determined by
the condition $[p^*g]=[\tilde{g}]$, where $p$ denotes the covering map
$\tM\rightarrow M$.

Assume now that $\Gamma\subset\mathcal{H}(\tM)$. Then
 the induced structure is globally conformal K\"ahler if and only
 if $\rho(\Gamma)=1$.
\end{pro}

\begin{D} The action of $\Gamma$ can be seen as satisfying the
 hypothesis of the reduction theorem, so the first claim
follows. However we give a straightforward construction.

Let $\tilde{g}$ be one of the K\"ahler metrics of the structure of \tM.
Given
an atlas $\{U_i\}$
for the covering map $p$,
induce
a local K\"ahler metric $g_i$ on any $U_i$  by projecting $\tilde{g}$ restricted
to one of the connected components of $p^{-1}(U_i)$.
Then $g_i$ and $g_{j}$
differ by a conformal map on $U_i\cap U_{j}$,
hence by a partition of unity of $\{U_i\}$
one
can glue the set $\{g_i\}$
to  a global metric $g$ which is
locally conformal the $g_i$'s, see Lemma~\ref{incollare},
hence is locally
conformal K\"ahler. The conformal class of
$g$ is uniquely defined by this construction. Moreover,
$p^*g$ is conformal
to $\tilde{g}$, as they are conformal on each component of the
covering $\{p^*(U_i)\}$ and again Lemma~\ref{incollare} holds.
If $g'$ is a Hermitian metric on $M$ such that $p^*g'$ is conformal to
$\tilde{g}$, then on each $U_i$ the restricted metric $g'|_{U_i}$ is conformal to
$g_i$ hence to $g|_{U_i}$, so $g$ and $g'$ are conformal,
again see Lemma~\ref{incollare}.

Now assume that $\Gamma\subset\mathcal{H}(K)$, and that
$\rho(\Gamma)\neq1$. Then $\Gamma$ is {\em not} contained in the
isometries of any K\"ahler metric of \tM. If in the class of $[g]$ there
existed
a K\"ahler metric $\bar{g}$ then its pull-back $p^*\bar{g}$ would belong to
$\langle \tilde{g}\rangle$. But $p^*\bar{g}$ being a pull-back implies that
$\Gamma$ acts with isometries with respect to it, which is absurd
since $\rho(\Gamma)\neq1$. Conversely, if $\rho(\Gamma)=1$ then $p$ is
a Riemannian covering space and $g$ itself is K\"ahler.
Hence the
induced locally conformal K\"ahler structure is globally conformal
K\"ahler if and only if $\rho(\Gamma)=1$.
\end{D}

This allows, under a natural condition, to compute locally
conformal K\"ahler reduction as having a
K\"ahler reduction as covering space. First remark that any group
$G\subseteq\aut{M}$ lifts to subgroups $\tG\subseteq\mathcal{H}(\tM)$ all having the
property that $p\circ \tG=G$.

\begin{teo}\label{compatibilitakaehler}
Let $(M,[g],J)$ be a locally conformal K\"ahler manifold, let
$G\subset\aut{M}$ be a subgroup satisfying the hypothesis of the
reduction theorem, and
admitting a lifting $\tG$ such that $\rho(\tG)=1$.
Then the K\"ahler reduction is defined, with momentum
map  denoted by $\moment_{\tM}$,  $\tG$ commutes
with the action of $\pi_1(M)$, and the following equality of locally conformal
K\"ahler structures holds:
\begin{equation}\label{ugualikaehler}
\moment^{-1}(0)/G\iso (\moment_{\tM}^{-1}(0)/\tG)/\pi_1(M).
\end{equation}

Conversely, let $\tG$ be a subgroup of isometries of a homothetic
K\"ahler manifold $(\tM,\langle \tilde{g}\rangle,J)$ of complex
dimension bigger than 1 satisfying the hypothesis
of K\"ahler reduction
and commuting with the
action of a subgroup $\Gamma\subset\mathcal{H}(\tM)$ acting freely and
properly discontinuously and such that $\rho(\Gamma)\neq1$. Moreover,
assume that $\Gamma$ acts freely and properly discontinuously on
$\moment_{\tM}^{-1}(0)/G$.
Then $\tG$
induces a subgroup $G$ of $\aut{M}$,
$M$ being  the locally conformal K\"ahler manifold $\tM/\Gamma$, which
satisfies the hypothesis of the reduction theorem, and the
isomorphism~\myref{ugualikaehler} holds.
\end{teo}

\begin{D}
To show that K\"ahler reduction is defined, one has to show that
the action of $\tG$ is Hamiltonian with respect to the globally
conformal
K\"ahler structure of \tM.  First remark that the Lie algebra of \tG\ coincides with
that of $G$, that we denote by $\mathfrak{g}$ as usual, and that the
fundamental vector field associated with $X\in\mathfrak{g}$ on \tM\ is
$p^*X$, where, as we claimed, we identify  $X$ with its associated
fundamental field on $M$.
Fix a metric $g\in[g]$ with Lee form $\lee$ and fundamental form
\fund\ and let $\moment$ be
the momentum map for $g$.
% Let $\tau\in C^\infty(\tM)$ be such that
%$d\tau=p^*\lee$.
Then we claim that
\dismap{\moment_{\tM}^\cdot}{\mathfrak{g}}{C^\infty(\tM)}{X}{p^*\moment^X}
is a momentum map for the action of \tG\ on \tM\ with respect to the
globally conformal K\"ahler metric $p^*g$. Indeed
\[
\begin{split}
d^{p^*\lee}\moment_{\tM}^{X}&=d^{p^*\lee}(p^*\moment^X)\\
&=p^*d^{\lee}\moment^X\\
&=p^*\iota_X\fund\\
&=\iota_{p^*X}p^*\fund.
\end{split}
\]
The same way one shows that $\moment_{\tM}$ is a homomorphism of
Poisson algebras, since such is $\moment_M$.
But now recall that from Remark~\ref{conformalpoisson} the property of
an action
to be twisted Hamiltonian is a conformal one, so the action of $G$ is
also twisted Hamiltonian for the K\"ahler metrics conformal to $p^*g$,
and is then ordinarily Hamiltonian for these K\"ahler metrics from
Proposition~\ref{globally}. This in turn implies, since $\rho(G)=1$,
that K\"ahler reduction is defined and
$\moment^{-1}(0)$ is diffeomorphic to $\moment_{\tM}^{-1}(0)/\pi_1(M)$.

As the action of $\tG$ is induced by $p$, it commutes with
the action of $\pi_1(M)$, so  the
following diagram of differentiable manifolds commutes:
\[
\xymatrix{
\moment^{-1}_{\tM}\ar[r]\ar[d]^p
&\moment^{-1}_{\tM}(0)/\tG    \ar[d]\\
\moment^{-1}(0)\ar[r]
&(\moment^{-1}_{\tM}(0)/\tG)/\pi_1(M)\iso\moment^{-1}(0)/G.\\
}
\]
Moreover the locally conformal K\"ahler structures induced on
$\moment^{-1}(0)/G$, as covered by the K\"ahler reduction $\moment_{\tM}^{-1}(0)/G$ and
as locally conformal K\"ahler reduction, are easily seen to coincide, and
 this ends the first part of the proof.

Conversely note that, as in Remark~\ref{coveringhomothetic},
if $\tau$ is such that $p^*\lee=d\tau$, then
$e^{-\tau} p^*{g}$ is K\"ahler, hence conformal to $\tilde{g}$. So $\Gamma$
acts as isometries of $e^{\tau}\tilde{g}$.
We claim that $e^\tau\moment_{\tM}$ is $\Gamma$-invariant, where by
$\moment_{\tM}$ we denote the K\"ahler momentum map. Postponing for
the moment the proof, this defines the locally
conformal K\"ahler momentum map as
$\moment^X\ug(e^\tau\moment_{\tM}^X)/\Gamma$, where we identify the Lee
algebra of $G$ with that of $\tG$.
This induced momentum map is easily shown to be
a homomorphism with respect to the Poisson structure. Moreover $\moment_M^{-1}(0)\iso
\moment_{\tM}^{-1}(0)/\Gamma$. Finally since the
action of $\Gamma$ on $\moment_{\tM}^{-1}(0)/G$ is free and properly
discontinuous, $(\moment_{\tM}^{-1}(0)/G)/\Gamma$ is a manifold, and since
the diagram
\[
\xymatrix{
\moment_{\tM}^{-1}(0)\ar[r]\ar[d]^p
&\moment_{\tM}^{-1}(0)/G\ar[d]\\
\moment_M^{-1}(0)\ar[r]
&(\moment_{\tM}^{-1}(0)/G)/\Gamma=\moment_M^{-1}(0)/G
}
\]
commutes, the action of $G$ on $\moment_M^{-1}(0)$ is proper and
free, that is, $G$ satisfies the hypothesis of the reduction
theorem.
 So the first part of the theorem implies the second.

We are left to prove the claim. For simplicity we write $\moment$
instead of $\moment_{\tM}$.
First we show that for any
$\gamma\in\Gamma$ and $X\in\mathfrak{g}$ we have
$\gamma^*\moment^X=\rho(\gamma)\moment^X$.
Indeed, recalling that $\gamma^*X=X$
\[
\begin{split}
d\gamma^*\moment^X&=\gamma^*d\moment^X\\
&=\gamma^*\iota_X\tfund\\
&=\iota_X\gamma^*\tfund\\
&=\rho(\gamma)\iota_X\tfund\\
&=\rho(\gamma)d\moment^X,
\end{split}
\]
hence $\gamma^*\moment^X-\rho(\gamma)\moment^X$ is constant on
$M$, and is equal to 0 since it is so on $\moment^{-1}(0)$.
So
$\gamma^*(e^\tau\moment^X)=e^{\gamma^*\tau}\rho(\gamma)\moment^X$.
But now recall that, from one side, the formula
$\gamma^*e^\tau\tfund=e^{\gamma^*\tau}\rho(\gamma)\tfund$ holds true,
from the other that
$\Gamma$ acts as isometries of $e^\tau\tilde{g}$, hence
$e^{\gamma^*\tau}\rho(\gamma)=1$, so the claim is true.
\end{D}

\section{Reduction of compact Vaisman manifolds}\label{vaisman}

\subsection{A conformal definition of compact Vaisman manifolds}

The original definition of Vaisman manifold is relative to a
Hermitian
manifold: the metric of a Hermitian manifold $(M,g,J)$ is a  {\em Vaisman
metric}\/ if it is locally conformal K\"ahler with $\omega$ non-exact
and
if $\nabla^g\lee_g=0$, where $\nabla^g$ is the Levi-Civita
connection.

\begin{defi} A conformal Hermitian manifold $(M,[g],J)$ is a
{\em Vaisman manifold}\/ if it is a locally conformal K\"ahler manifold,
non globally conformal K\"ahler and admitting a Vaisman metric in
$[g]$.
\end{defi}

The condition on  the parallelism of the Lee form
 is {\em not}\/
invariant up to
conformal changes of metric, and there is not in the literature a
conformally invariant
criterion to decide whether a given locally conformal K\"ahler
manifold is Vaisman.

Such a  criterion was recently given in \cite{KaOGFC} in the case of compact
locally conformal K\"ahler manifolds. Here we shall use it to derive a
presentation for compact Vaisman manifolds that behaves effectively with respect to reduction.

The
construction strictly links Vaisman geometry with Sasaki
geometry. We start with the following definition-proposition, which is
equivalent to the standard one. On this subject see~\cite{BlaRGC,BoGTSM}.

\begin{defi}
Let $(W,g_W,\eta)$ be a Riemannian manifold of odd dimension bigger
than 1
with a contact form $\eta$ such that on the distribution $\eta=0$
the $(1,1)$-tensor
$J$ that associates to a vector field $V$ the vector field
$\sharp_g\iota_V d\eta$ satisfies $J^2=-1$. Call $\zeta$ the Reeb
vector field of $\eta$.
Define on  the cone $W\times\R$
the metric $g=e^{t}dt\otimes dt+e^{t}\pi^*g_W$ and the complex structure
that extends $J$ associating to $dt$ the vector field
$\pi^*\zeta$, $\pi$ being the projection of $W\times\R$ to $W$. This
is equivalent to assigning to $W\times\R$ the same $J$ and the
compatible symplectic form $\fund\ug d(e^t\pi^*\eta)$.
Then we say that $(W,g_W,\eta)$ is a {\em Sasaki manifold} if its cone
$(W\times \R,e^{t}dt\otimes dt+e^{t}\pi^*g_W,J)$ is
K\"ahler.
\end{defi}

The standard example is that of the odd-dimensional sphere
contained in $\C^n\setminus 0$, with $n\geq2$. The usual
K\"ahler metric $\sum dz_i\otimes d\bar{z}_i$ associated to
the complex structure of $\C^n\setminus 0$ restricts to the sphere to a Riemannian
metric and a CR-structure, respectively, that give to
$S^{2n-1}$ the Sasaki structure whose cone is
$\C^n\setminus 0$ itself, via the identification $(x,t)\mapsto
e^{t/2}x$.

It is well-known that the conformal metric $|z|^{-2}\sum dz_i\otimes
d\bar{z}_i$ has parallel Lee form. This property extends to every
K\"ahler cone, as is implicit in \cite{KaOGFC}.

\begin{lem}\label{parallela}
The K\"ahler cone $(W\times\R,g,J)$ of a Sasaki manifold
admits the metric $\tilde{g}=2e^{-t}g$ in its
conformal class such that $\nabla^{\tilde{g}}\lee_{\tilde{g}}=0$.
In particular the Lee form of $\tilde{g}$ is $-dt$.

\end{lem}

\begin{D}
Recall that the fundamental form \fund\ of $g$ is K\"ahler, so the
fundamental form $\tfund= 2e^{-t}\fund$ of $\tilde{g}$ is such that
$d\tfund=-2e^{-t}dt\wedge\fund=-dt\wedge\fund$. So $-dt$ is the
Lee form of $\tilde{g}$. Remark that
\[
\tilde{g}=2dt\otimes dt + 2\pi^*g_W
\]
This shows that $\partial_t$ is twice the metric dual of $-dt$.
Recall that for any 1-form $\sigma$ the following holds
\[
2\tilde{g}(\nabla^{\tilde{g}}_X\sigma^\sharp,Z)=(\lie_{\sigma^\sharp}\tilde{g})(X,Z)+d\sigma
(X,Z)
\]
where $\lie$ denotes the Lie derivative. So, in our case, since
$dt$ is closed, we only have to show that
$\partial_t$ is Killing. But this is true since
$\lie_{\partial_t}\tilde{g}=2\lie_{\partial_t}dt\otimes dt$, and
$\lie_{\partial_t}dt=d\iota_{\partial_t}dt=0$.
\end{D}

The following also follows from computations developed in \cite{KaOGFC}.

\begin{pro}\label{quoziente}
Let $(W,g_W,\eta)$ be a Sasaki manifold, let $\Gamma$ be a subgroup of
$\mathcal{H}(W\times\R)$ acting freely and properly
discontinuously on $W\times\R$, in such a way that
$\rho(\Gamma)\neq1$ and for any $\gamma\in\Gamma$
\[
\gamma\circ\phi_t=\phi_t\circ\gamma,
\]
that is, $\Gamma$ commutes with the real flow generated by $\partial_t$.

Then the induced locally conformal K\"ahler structure on
 $M\ug(W\times\R)/\Gamma$ is a Vaisman structure,
not globally conformal K\"ahler.
\end{pro}

\begin{D}
Since $\Gamma\subset\mathcal{H}(W\times\R)$ and $\rho(\Gamma)\neq1$
the quotient has a locally non globally
conformal K\"ahler structure, recall Proposition~\ref{indurre}.

To show that the structure is Vaisman we show that $\Gamma$ acts
by isometries of the metric $\tilde{g}=2e^{-t}g_{W\times\R}$, where by
$g_{W\times\R}$ we denote the cone metric on $W\times\R$.
% From Remark~\ref{coveringhomothetic} this is equivalent to show that
% $p^*\bar{\omega}=-dt$, where by $\bar{\lee}$ we denote the Lee form of
% the quotient metric, since
% \[
% e^tp^*g=\lambda g_{W\times\R}
% \]
% that is, $\Gamma$ act as isometries of $e^{-t}g_{W\times\R$.
%
This is
equivalent to show that $\Gamma$ acts by symplectomorphisms of the
conformal K\"ahler form $2e^{-t}d(e^t\pi^*\eta)$.

We claim that for any $\gamma\in\Gamma$ the following properties hold:
\[
\begin{split}
\gamma^*\pi^*\eta=\pi^*\eta\\
\gamma^*e^t=\rho(\gamma)e^t.
\end{split}
\]
For this, first note that $\gamma$ commuting with the real natural
flow implies $\gamma_*\partial_t=\partial_t$, it being holomorphic
implies $\gamma_*J\partial_t=J\partial_t$ and it being conformal
implies $\gamma_*\langle\partial_t,J\partial_t\rangle^\perp
=\langle\partial_t,J\partial_t\rangle^\perp$.
Now remark that for $X\in\langle\partial_t,J\partial_t\rangle^\perp$
$\pi_*X\in \nul\eta$, so
\[
\gamma^*\pi^*\eta(X)=\eta(\pi_*\gamma_*X)=0.
\]
Moreover $1=\eta(\zeta)=\eta(\pi_*(J\partial_t))=\eta(\pi_*(\gamma_*J\partial_t)=\gamma^*\pi^*\eta(J\partial_t)$
and this implies the first claim. Now recall that
$\pi^*\eta=e^{-t}\iota_{\partial_t}\fund$, so
\[
\begin{split}
\gamma^*(e^t)\rho(\gamma)\iota_{\partial_t}\fund
&=\gamma^*(e^{-t}\iota_{\partial_t}\fund)\\
&=\gamma^*\pi^*\eta\\
&=\pi*\eta\\
&=e^{-t}\iota_{\partial_t}\fund
\end{split}
\]
which shows the second claim.

Then it follows
\[
\begin{split}
\gamma^*(2e^{-t}d(e^t\pi^*\eta))&=2\rho(\gamma)^{-1}e^{-t}d\gamma^*(e^t\pi^*\eta)\\
&=2\rho(\gamma)^{-1}e^{-t}\rho(\gamma)d(e^t\gamma^*\pi^*\eta)\\
&=2e^{-t}d(e^t\pi^*\eta).
\end{split}
\]

So $\tilde{g}$ factors through the action of $\Gamma$, hence
inducing $g_M$ on $M$ which, by Lemma~\ref{parallela}, is Vaisman, and belongs to the  locally conformal K\"ahler
structure of $M$ since $p^*g_M=\tilde{g}\sim g_{W\times\R}$.
\end{D}

The characterization  given in
\cite{KaOGFC} shows in fact that {\em
any} compact Vaisman manifold is produced this way. We briefly recall
this construction, since some details which are less relevant in that
work become necessary in this one, so we need to express them explicitly.

Remark that a vector field $V$ generating  a 1-parameter subgroup of
\aut{M} does {\em not} imply that the flow of $JV$ is contained in
\aut{M}. If this happens, the set of the flows of the subalgebra
generated by $V$ and $JV$ is a Lie subgroup of \aut{M} of real
dimension 2 that has a structure of complex Lie group of
dimension 1.
This motivates the following definition:

\begin{defi}[\cite{KaOGFC}]~
A {\em holomorphic conformal flow} on a locally conformal
K\"ahler manifold $(M,[g],J)$ is a 1-dimensional complex Lie subgroup of
the biholomorphisms of $(M,J)$ which is contained in $\aut{M}$.
\end{defi}

\begin{oss}
The field $\partial_t$ on a K\"ahler cone of a Sasaki manifold
generates a holomorphic conformal flow. Its flow
$\phi_s(w,t)=(w,t+s)$ is in fact contained
in $\mathcal{H}(W\times\R)$, and satisfies $\rho(\phi_s)=e^s$, since
$\pi\phi_s=\pi$ and
\[
\begin{split}
\phi^*_s(d(e^t\pi^*\eta))&=d(e^{t+s}\phi_s^*\pi^*\eta)\\
&=e^sd(e^t\pi^*\eta).
\end{split}
\]
 The
flow of $J\partial_t$, which is a vector field that restricts to the Reeb vector
field of $W$, which is  a Killing vector field of $W$, generates
 isometries of $W\times\R$.
We call the real flow generated by $\partial_t$ the {\em natural real flow}
and the holomorphic conformal flow generated by $\partial_t$ the
{\em natural holomorphic flow} of the K\"ahler cone.

Finally remark that for a biholomorphism $h$ of a Hermitian manifold
to commute with
the flow of a vector field $V$ it is necessary and sufficient that it
commutes with the whole
holomorphic flow, since $h^*V=V$ is equivalent to $h^*JV=JV$. So if a
holomorphic conformal flow $\mathcal{C}$ is defined on a locally conformal
K\"ahler manifold saying that it is preserved by an automorphism is $h$
equivalent to saying that $h$ preserves a real generator
of $\mathcal{C}$.
\end{oss}

\begin{teo}[Kamishima \& Ornea, \cite{KaOGFC}]\label{characterizationvaisman}
Let $(M,[g],J)$ be a compact, connected, non-K\"ahler, locally
conformal K\"ahler manifold of complex dimension $n\geq2$.
Then $(M,[g],J)$ is Vaisman if and only if $\aut{M}$
admits a holomorphic conformal flow.
\end{teo}
\begin{D} For the technical lemmas we refer directly to the cited paper.

First, if $M$ is a Vaisman manifold, then the dual vector field $\sharp\lee$ of
its Vaisman metric generates a holomorphic conformal flow, that
is, both the flow of $\lee^\sharp$ and the flow of $J\lee^\sharp$
belong to \aut{M}, see~\cite{DrOLCK}.

On the opposite direction let $\mathcal{C}$ be the holomorphic
conformal flow on $M$. Fix a lift $\tilde{\mathcal{C}}$ of
$\mathcal{C}$ to $\tilde{M}$. One proves (Lemma 2.1) that $\rho(\tC)=\R_+$.
Choose a vector field $\xi$ on \tM\ such that the flow $\{\psi_t\}$ of $-J\xi$ is
contained in $\ker\rho|_{\tC}$.
Remark that the flow $\{\phi_t\}$ of $\xi$
is also contained in \tC, since \tC\ is a holomorphic conformal
flow, and that
$t\mapsto\rho(\phi_t))$ is surjective. Choose
$\tfund$ in the homothety class of K\"ahler forms on \tM\ in such
a way that this homomorphism is $t\mapsto e^t$, that is for any
$t$
\[
\phi_t^*\tfund=e^t\tfund\qquad \psi_t^*\tfund=\tfund.
\]
In particular the subgroup $\{\phi_t\}$ of \tC\ is isomorphic to
\R.

At this step the hypothesis of compactness of $M$ is crucial:
using this fact
one proves that the action of $\{\phi_t\}$ is free and proper (Lemma 2.2). In
particular $\xi$ is never vanishing.

Define the smooth map
\dismap{s}{\tM}{\R}{x}{\tfund(J\xi_x,\xi_x)}
and remark that $1$ is a regular value of $s$, that $s^{-1}(1)$ is
non empty, hence is a regular submanifold of \tM\ that we denote
by $W$ (Proposition 2.2). Note that $W$ is the submanifolds of
those points where $\xi$ has unitary norm. In particular one
proves that if $x$ is in $W$ then $d_xs(\xi_x)=1$, so $\xi$ is
transversal to $W$.

Denote by $i$ the inclusion of $W$ in $M$.
It turns then out that $(W,i^*\iota_\xi\tfund, i^*g)$ is a connected Sasaki
manifold,
and
that
\dismap{H}{W\times\R}{\tM}{(w,t)}{\phi_t(w)}
is an isometry with respect to the K\"ahler cone structure on
$W\times\R$.

One is left to show that $\pi_1(M)$ satisfies the conditions of
Proposition~\ref{quoziente}. Indeed $\rho(\pi_1(M))\neq1$ since $M$ is
non K\"ahler, and $\pi_1(M)$ commutes with the real flow generated by
$\partial_t$ since this last factors to $M$.
\end{D}

The proof of this theorem proves in particular the following fact.

\begin{cor}
Any compact Vaisman manifold $(M,[g],J)$ can be presented as
a pair $(W,\Gamma)$ where $W$ is a compact Sasaki manifold and
$\Gamma\subset\mathcal{H}(W\times\R)$ such that $M$ is isomorphic as a
locally conformal K\"ahler manifold to $(W\times\R)/\Gamma$. Moreover
$W$ can be chosen to be simply connected, hence $W\times\R$ is the
universal covering of $M$ and $\Gamma$ is isomorphic to $\pi_1(M)$.
\end{cor}

\begin{oss}
This can be reformulated in the following way. Consider the collection
of pairs
$(W,\Gamma)$ as in the previous corollary. Given a Sasaki morphism
one naturally induces a morphism  on the K\"ahler cones by composing with identity on the
factor \R. So define the category $\mathcal{S}$ of pairs $(W,\Gamma)$ by considering
as morphisms between $(W,\Gamma)$ and $(W',\Gamma')$ those Sasaki
morphisms (i.e. isometries preserving the contact form)
which induce between $W\times\R$ and $W'\times\R$ morphisms
which are equivariant with respect to the actions of $\Gamma$ and
$\Gamma'$. Define a functor between this category and the category of
Vaisman manifolds (with morphisms given by holomorphic, conformal
maps) by associating to $(W,\Gamma)$ the manifold $(W\times\R)/\Gamma$
and to a morphism the induced morphism between the quotients.
This functor is surjective on the objects of the
subcategory of compact Vaisman manifolds, but not on the morphisms.
\end{oss}

We say that $(W,\Gamma)$ is a {\em presentation} of $M$ if $M$ is in
the image of $(W,\Gamma)$ by this functor. %Conversely given $M$ compact Vaisman manifold the
%{\em associated Sasaki manifold} is the one found in the proof of
%Theorem~\ref{characterizationvaisman} by means of the {\em Vaisman
%holomorphic flow} of $M$, that is, the holomorphic flow generated by the
%Lee field.
This construction suggests the following definition.

\begin{defi} Let $(M,[g],J)$ be a compact Vaisman manifold. A locally conformal K\"ahler morphism is
a {\em Vaisman morphism} if it belongs to the image of the functor from $\mathcal{S}$ to Vaisman manifolds.
In particular {\em Vaisman automorphisms} are a subgroup of
\aut{M,[g],J}.
\end{defi}

Rermark that Vaisman morphisms are
the set of locally conformal K\"ahler  morphisms $h$ that commute with the {\em Vaisman
real flow}, that is, the real flow generated by the Lee field, and
that admit a lifting
 $\tilde{h}\in\mathcal{H}(\tM)$ such that $\rho(\tilde{h})=1$.

\begin{oss}\label{isometrie} It must
be noted that whenever $(M,[g],J)$ is Vaisman and compact the Gauduchon metric is the
Vaisman metric (well-defined up to homothety). This implies that
\aut{M}
{\em coincides} with the holomorphic isometries of the Vaisman metric
in this case. So Vaisman automorphisms coincide with the set of
those $h$
being isometries of the Vaisman metric that commute with the Vaisman
real flow and admitting a lifting $\tilde{h}\in\mathcal{H}(\tM)$ such that $\rho(\tilde{h})=1$.
%
%A general
%isometry lifts to a homothety of $\tM\iso W\times\R$ that commutes with
%the natural flow.
%In \cite{KaOGFC}, Remark 2.1, it is shown that the subgroup
%$\mathcal{C}_{\mathcal{H}}$ of all the homotheties of $W\times\R$ that
%commute with the natural real flow splits by means of $\rho$ as
%the group
%$Isom(W)\times\R$, where $\R$ is the real flow itself.
%So if $M$ is compact Vaisman the elements of $\aut{M}$ lift to the
%subgroup of this product of those maps
%also commuting with $\pi_1(M)$.
\end{oss}

\begin{oss}  In passing we remark that the categories we are talking about all
admit a forgetting functor in corresponding ``symplectic''
categories, respectively locally conformal symplectic manifolds,
homothetic symplectic manifolds and contact manifolds (with the
functor given by symplectic cone construction). By associating to
the pairs $(W,\Gamma)$ the locally conformal symplectic manifolds
$(W\times\R)/\Gamma$ one obtains a subcategory that might represent
the symplectic version of Vaisman manifolds. See~\cite{VaiLCS} as
a leading reference on locally conformal symplectic manifolds, in
particular see the notion of locally conformal symplectic
manifolds of first kind.

\end{oss}

\subsection{Reduction for compact Vaisman manifolds}

It is noted in \cite{BiGCML} that $G$ acting by isometries with respect to a Vaisman
metric $g$ does
not imply that the reduced metric is Vaisman, since $\lee_g$ being
parallel with respect to the Levi-Civita connection of $g$
does not imply its restriction to $\moment^{-1}(0)$ being
parallel.

We prove that our reduction is compatible with Sasaki reduction,
see~\cite{GrORSM}, and thus show in Theorem~\ref{quozientevaisman} that reduction of compact Vaisman
manifolds
by the action of Vaisman automorphisms (whose action results to be always
twisted Hamiltonian)
produces Vaisman manifolds.
Given a Sasaki manifold $(W,\eta,J)$ we call {\em Sasaki isomorphisms}\/ and denote by $\isom{W}$ the
diffeomorphisms of $W$ preserving both the metric and the contact
form.

\begin{teo}\label{compatibilitasasaki}
Let $((W,\eta,J),\Gamma)$
be  a pair in the category $\mathcal{S}$ and denote by $M$ the
associated Vaisman manifold.
 Let
$G\subset \isom{W}$ be a subgroup  satisfying the hypothesis of Sasaki
reduction. Then $G$ can be considered as a subgroup of
$\mathcal{H}(W\times\R)$. Assume that the action of $G$ commutes with that of
$\Gamma$, and that $\Gamma$ acts freely and properly discontinuously
on the K\"ahler cone $(\moment_W^{-1}(0)/G)\times\R$.

Then
$G$ induces a subgroup of $\aut{M}$ satisfying the hypothesis of the
reduction theorem, and the reduced manifold is isomorphic with
$((\moment_W^{-1}(0)/G)\times\R)/\Gamma$. In particular the reduced
manifold is Vaisman.

\end{teo}

\begin{D}
We first prove that the induced action satisfies the hypothesis of the
K\"ahler reduction. The momentum map $\moment_W$ is
\dismap{\moment_W}{\mathfrak{g}}{C^\infty(W)}{X}{\iota_X\eta.}
Remark that the fundamental field associated to $X\in\mathfrak{g}$ on
$W\times\R$ is projectable, so we can define
$\moment_{W\times\R}^X\ug e^t\pi^*\iota_{\pi_*X}\eta$. To show that
this  is a momentum
map for the action of $G$ on $W\times\R$ we directly compute
\[
\begin{split}
d\moment_{W\times\R}^X&=d(e^t\pi^*\iota_{\pi_*X}\eta)\\
&=d(e^t\iota_X\pi^*\eta)\\
&=d(\iota_Xe^t\pi^*\eta)\\
&=\iota_Xd(e^t\pi^*\eta)-\lie_X(e^t\pi^*\eta)\\
&=\iota_X\fund
\end{split}
\]
where $\lie_X(e^t\pi^*\eta)=0$ comes from the properties of the
action. So the action of $G$ is weakly Hamiltonian. Moreover
\[
\begin{split}
\{\moment_{W\times\R}^X,\moment_{W\times\R}^Y\}&=\fund(X,Y))\\
&=d(e^t\pi^*\eta)(X,Y)\\
&=Y(e^t\pi^*\eta(X))-X(e^t\pi^*\eta(Y))+\iota_{[X,Y]}e^t\pi^*\eta\\
&=0+0+e^t\pi^*\iota{\pi_*[X,Y]}\eta\\
&=\moment_{W\times\R}^{[X,Y]}
\end{split}
\]
where $Y(e^t\pi^*\eta(X))=X(e^t\pi^*\eta(Y))=0$ is due to the
properties of the action. The action of $G$ is then Hamiltonian, and
K\"ahler reduction is defined.

So
$\moment_{W\times\R}^{-1}(0)\iso(\moment_W^{-1}(0)\times\R)$, and the
action of $G$ being proper and free on $\moment_W^{-1}(0)$ implies it
having the same properties on $\moment_{W\times\R}^{-1}(0)$. Moreover
the action of $\Gamma$ being free and properly discontinuous on
$(\moment_W^{-1}(0)/G)\times\R=\moment_{W\times\R}^{-1}(0)/G$ implies one
can apply Theorem~\ref{compatibilitakaehler}, so the isomorphism
is proven. Then by applying Proposition~\ref{quoziente} one sees
that the reduced manifold is Vaisman.
The theorem then follows.
\end{D}

\begin{oss}
The previous theorem also applies to non-compact Vaisman manifolds of
the form $(W\times\R)/\Gamma$.

\end{oss}

But now recall that Sasaki reduction, as contact reduction in
fact, does {\em not} need a notion of Hamiltonian action, that is,
the action of any subroup of Sasaki isometries lead to a momentum
map, hence to reduction, up to topological conditions, that is $\moment^{-1}(0)$ being
non-empty, 0 being a regular value for $\moment$, and $G$ acting
freely and properly discontinuously on $\moment^{-1}(0)$. This
implies the main result of this section.

\begin{teo}\label{quozientevaisman}
Let $(M,[g],J)$ be a compact Vaisman manifold. Let
$G\subset\aut{M}$ be a subgroup  of Vaisman automorphisms.
%Denote by $W$ the
%compact Sasaki manifold associated with $M$.
%
%Then $G$ induces a subgroup $\tG$ of
%isometries on $W$ and the action of $G$ on $M$ is twisted
%Hamiltonian. If $G$ acts freely and properly on $\moment^{-1}(0)$
%then
% it satisfies the hypothesis of Sasaki reduction
%and the locally conformal K\"ahler reduction is isomorphic to
%$((\moment_W^{-1}(0)/\tG)\times\R)/\pi_1(M)$. In particular the reduced
%manifold is Vaisman.
%
Then the action of $G$ on $M$ is twisted
Hamiltonian. If $\moment^{-1}(0)$ is non-empty, $0$ is a regular
value for $\moment$ and
$G$ acts freely and properly on $\moment^{-1}(0)$, then the
reduced manifold is Vaisman.

Moreover for any $(W,\Gamma)$  presentation of $M$ and for any $\tG$
subgroup of isometries of $W$ in
the preimage of $G$, the group $\tG$ induces a Sasaki reduction
$\moment_W^{-1}(0)/\tG$ and
the Vaisman reduced manifold is isomorphic to $((\moment_W^{-1}(0)/\tG)\times\R)/\Gamma$.
\end{teo}

\begin{D}
The automorphisms
of $M$ lift to a subgroup \tG\ of $\mathcal{H}(W\times\R)$ of maps commuting
with the natural flow and such that $\rho(\tG)=1$.
By definition of Vaisman automorphisms
$\tG$ is contained in $\isom{W}$. So the compatibility between
Sasaki and K\"ahler reduction of~\cite{GrORSM} applies, hence
$\moment_{W\times\R}^{-1}(0)/\tG$ is isomorphic with the K\"ahler cone $(\moment_W^{-1}(0)/\tG)\times\R$.
%\[
%\xymatrix{
%\moment_W^{-1}(0)\ar@{(_>}[rr]^i\ar[d]&&\moment_{\tM}^{-1}(0)\ar[ld]\ar[dd]_p\\
%\moment_W^{-1}(0)/G\ar[r]&\moment_{\tM}^{-1}(0)/G\ar[d]\\
%&\moment_M^{-1}(0)/G\ar[r]&\moment_M^{-1}(0)
%}
%\]

Moreover $\Gamma$ acts freely and properly
discontinuously on it, since the quotient $\moment_M^{-1}(0)/G$ is a manifold, and commutes with
its natural  real flow. Then Theorem~\ref{compatibilitasasaki}
applies,
and this proves
the theorem.
\end{D}

%\begin{oss}
%We stress that we showed that the property, typical of Sasaki (in
%fact, of contact) geometry of {\em always}\/ admitting a momentum
%map is conserved under the functor between the category
%$\mathcal{S}$ and the category of Vaisman manifolds.
%\end{oss}

\section{A class of examples: weighted actions on Hopf
manifolds}\label{vai}

We apply the theorem in last section to the simple case when $\Gamma=\Z$ is contained as a discrete subgroup of the
natural holomorphic flow of $W\times\R$. In particular the Vaisman
manifold topologically
is simply $W\times S^1$.
This nevertheless covers much of the already known examples of Vaisman
manifolds, as was shown in \cite{KaOGFC}.
First consider $S^{2n-1}$  equipped with the standard CR structure
$J$
coming from
$\C^n$.
The action on $\C^n\setminus0$ of  the cyclic group $\Gamma_\alpha$ generated
by $z\mapsto\alpha z$  (for any $\alpha\in\C$ such that $|\alpha|>1$)
produces the so-called {\em standard Hopf manifolds}. For any $\{c_1,\dots,c_n\}\in (S^1)^n$
and any  set $A\ug\{a_1,\dots,a_n\}$ of real numbers such that
$0<a_1\leq\dots\leq a_n$
the action of the cyclic group
$\Gamma_{\{c_1,\dots,c_n\},A}\subset\mathcal{C}$
generated by $(z_1\dots z_n)\mapsto
(e^{a_1}c_1z_1,\dots,e^{a_n}c_nz_n)$ produces the complex
manifolds  usually called  {\em non-standard Hopf manifold}.

Let $\eta_0$ be the
Sasaki
structure
coming from the standard form $\fund=-i\sum dz_i\wedge d\bar{z}_i$
of $\C^n$. The action of any $\Gamma_\alpha$ is
is by homotheties for the cone structure, hence produces Vaisman
structures $((S^{2n-1},\eta_0,J),\Gamma_\alpha)$ on standard Hopf manifolds.
More generally, for any $A\ug\{a_1,\dots,a_n\}$ of real numbers such that
$0<a_1\leq\dots\leq a_n$ let $\eta_A$ be defined  the following way:
\[
\eta_A\ug\frac{1}{\sum a_i|z_i|^2}\eta_0.
\]
Fixed $A$ one obtains that for any $\{c_1,\dots,c_n\}\in (S^1)^n$ the action of $\Gamma_{\{c_1,\dots,c_n\},A}$
 is by homotheties of the
corresponding cone structure on $\C^n\setminus0$, hence inducing
Vaisman structures
\[
((S^{2n-1},\eta_A,J),\Gamma_{\{c_1,\dots,c_n\},A})
\]
 on the
non-standard Hopf manifolds
 (cf.\ \cite{KaOGFC}).

So if we act on $(S^{2n-1},\eta_A,J)$ by a circle of Sasaki isometries and $n>2$ we
generate a Vaisman reduced manifold of dimension $2n-2$ for every
$\Gamma_{\{c_1,\dots,c_n\},A}$, whose underlying manifold is the
product of the
Sasaki reduced manifold with $S^1$.

Remark that the contact structures of the Sasaki
manifolds $(S^{2n-1},\eta_A,J)$ all coincide. Denote by
$\cont{S^{2n-1}}$ the set of contact automorphisms of $S^{2n-1}$, which
simply coincide with restriction of biholomorphisms of $\C^n$.

For any $\Lambda=(\lambda_1,\dots,\lambda_n)\in \R^n$ let
$G_{\Lambda}\subset \cont{S^{2n-1}}$ be the subgroup of
those maps $h_{\Lambda,t}$, $t\in\R$, such that
\[
h_{\Lambda,t}(z_1,\dots,z_n)=(e^{ i\lambda_1t}z_1,\dots,e^{ i\lambda_nt}z_n).
\]
Remark that any $G_{\Lambda}$ is composed in
fact by holomorphic isometries of the standard  K\"ahler
structure. Moreover a direct computation shows that its action on
$S^{2n-1}$
is by isometries for any of the $\eta_A$.
We
call the action of $G_\Lambda$
{\em weighted}\, by the {\em weights} $(\lambda_1,\dots,\lambda_n)$.
 We restrict to the
$\Lambda$'s such that
$G_{\Lambda}$ is isomorphic to $S^1$: it is easy to see that this
happens whenever
the ratios between the weights
are rational.

 The corresponding momentum map for the Sasaki manifold $(S^{2n-1},\eta_A,J)$ is defined by:
\[
\moment_{\Lambda}^1(z)=H_{\Lambda}(z)\ug\frac{1}{2(\sum a_i|z_i|^2)}(\lambda_1|z_1|^2+\dots+\lambda_n|z_n|^2).
\]
So a Sasaki reduction is defined whenever the weights are such
that $\moment^{-1}(0)$ is not empty and
the action on $\moment^{-1}(0)$ is free and proper. The condition that
$\moment^{-1}(0)$ is not empty is equivalent to requiring that the
signs of the $\lambda_i$ are not all the same.

 Let $k\in \{1,\dots,n-1\}$ be the number of negative
weights of $\Lambda$, and assume the negative weights are the first $k$.
Then there is a diffeomorphism
\dismap{\Phi_{\Lambda}}{
S^{2k-1}\times S^{2n-2k-1}}{\moment^{-1}(0)}{((\xi_1,\dots,\xi_k),(\zeta_1,\dots,\zeta_{n-k}))}
{(\frac{\xi_1}{\sqrt{-\lambda_1}},\dots,\frac{\xi_k}{\sqrt{-\lambda_k}},\frac{\zeta_1}{\sqrt{\lambda_{k+1}}},
\dots,\frac{\zeta_{n-k}}{\sqrt{\lambda_{n}}})}
equivariant with respect to the
action
$$
w_{\lambda,t}((\xi_1,\dots,\xi_k),(\zeta_1,\dots,\zeta_{n-k}))=
((e^{ i\lambda_1t}\xi_1,\dots,e^{ i\lambda_kt}\xi_k),(e^{ i\lambda_{k+1}t}\zeta_1,\dots,
e^{ i\lambda_nt}\zeta_{n-k}))
$$
from one side
and the action of $G_{\Lambda}$ on
$\moment^{-1}(0)$ from the other: $h_{\Lambda,t}\circ
\Phi_\Lambda=\Phi_\Lambda\circ w_{\lambda,t}$.

Call $S(\Lambda)$ the quotient of this action. This will generally be an orbifold. A sufficient
condition for the action of $G_\Lambda$  to be free, hence for
$S(\Lambda)$ to be a manifold, is that the $\lambda_i$'s are
relatively prime integers. Recall that for any
$\Lambda$ such that $S(\Lambda)$ is a manifold Theorem~\ref{compatibilitasasaki} implies there
exists a Vaisman structure on $S(\Lambda)\times S^1$ for every
$(a_1,\dots,a_n)\in\R^n$ such that $0<a_1\leq\dots\leq a_n$ and for every
$(c_1,\dots,c_n)\in(S^1)^n$, each being the reduction of the Hopf
manifold associated to $(a_1,\dots,a_n)$ and with Vaisman structure associated
to $((a_1,\dots,a_n),(c_1,\dots,c_n))$.

We analyze the topological type of
the reductions $S(\Lambda)$ in some of these cases.

\begin{exa} Assume that $n\geq2$, $k=1$, that is, $\lambda_1<0$, $\lambda_i>0, i=2,\dots,n$. Then the space
$\moment^{-1}(0)$ is diffeomorphic to $S^1\times
S^{2n-3}$.

One easily shows that $S(-1,1,\dots,1)$ is $S^{2n-3}$. One can also show that
the Sasaki structure reduced form the standard is again the standard
one, so any standard Hopf manifold comes as a reduction of the
corresponding Hopf manifold of higher dimension. This is also
shown in~\cite{BiGCML}.

In turn, as shown in~\cite{GrORSM}, for  any
negative integer $p$, $S(p,1,\ldots,1)$  is diffeomorphic to
$S^{2n-3}/\Z_p$, so we obtain a family of  Vaisman
structures on $(S^{2n-3}/\Z_p)\times S^1$. In particular for $n=3$
we obtain Sasaki structures on lens spaces of the form $L(p,1)$, hence Vaisman
structures on complex surfaces diffeomorphic with $L(p,1)\times S^1$.

\end{exa}

\begin{exa}
If $n=4$, $k=2$, then $\moment^{-1}(0)$ is diffeomorphic
to $S^3\times
S^{3}$.
In particular $S(-1,-1,1,1)$ is known to be $S^3\times S^2$, see
\cite{GrORSM}. So we obtain a family of Vaisman structures
on $S^3\times S^2 \times S^1$. It is
interesting to note that reducing from the standard structure one
obtains a
manifold that also bears a semi-K\"ahler
structure, when seen as twistor space of the standard Hopf surface.
\end{exa}

\begin{exa}
More generally if $n=2k=4s$, $\moment^{-1}(0)$ is
diffeomorphic to $S^{4s-1}\times S^{4s-1}$. In analogy with the
case $n=4$ treated in \cite{GrORSM}
apply the
diffeomorphism
\dismapnoname{S^{2k-1}\times S^{2k-1}}{S^{2k-1}\times S^{2k-1}}{((\xi_1,\dots,\xi_k),(\zeta_1,\dots,\zeta_{k}))}{\hspace{-9.25pt}((\xi_1\zeta_k+\overline{\xi_2\zeta_{k-1}},\xi_1\zeta_{k-1}-\overline{\xi_2\zeta_{k}},\dots,
\xi_{k-1}\zeta_2+\overline{\xi_k\zeta_{1}},\xi_{k-1}\zeta_{k}-\overline{\xi_k\zeta_{2}}),
(\zeta_1,\dots,\zeta_{k}))}
and remark that it is equivariant with respect to the action of
$G_{(-1,\dots,-1,1,\dots,1)}$ on the first space and the action on
the second space given by the product of
the trivial action on the first factor and the Hopf action on the second factor. This
proves that $S(-1,\dots,-1,1,\dots,1)$ is diffeomorphic to
$S^{4s-1}\times \Cp{2s-1}$. Thus we obtain families of Vaisman
structures on $S^{4s-1}\times \Cp{2s-1}\times S^1$.
%$S(-1,-1,-1,-1,1,1,1,1)$ is diffeomorphic to
%$S^6\times S^7$ and one obtains a family of Vaisman structures
%on $S^1\times S^6 \times S^7$.
%
%These families of Vaisman structures include the ones previously
%found in \cite{OrPIHB}.
\end{exa}

%\nocite{}
%\bibliography{rosa,mau}

%\end{document}

\begin{address}

Rosa Gini \& Maurizio Parton, Dipartimento di Matematica,
Universit\`a di Pisa, via Buonarroti $2$, I-56127 Pisa, Italy.
E-mail addresses: {\tt gini@dm.unipi.it {\rm \&} parton@dm.unipi.it}.

Liviu Ornea,
University of Bucharest,
Faculty of Mathematics,
70109 14 Academiei str.,
Bucharest, Romania. E-mail address:
{\tt lornea@imar.ro}.

\end{address}

\end{document}